\documentstyle[10pt,amscd,xypic,amssymb]{amsart}

\xyoption{all}
\CompileMatrices

\newcommand{\nc}{\newcommand}

\makeatletter
\@addtoreset{equation}{section}
\makeatother

\newtheorem{theorem}[equation]{Theorem}
\newtheorem{proposition}[equation]{Proposition}
\newtheorem{lemma}[equation]{Lemma}
\newtheorem{corollary}[equation]{Corollary}
\newtheorem{conjecture}[equation]{Conjecture}

\theoremstyle{definition}

\newtheorem{definition}[equation]{Definition}

\theoremstyle{remark}

\newtheorem{remark}[equation]{Remark}

\nc{\BA}{{\mathbb{A}}}
\nc{\BC}{{\mathbb{C}}}
\nc{\BM}{{\mathbb{M}}}
\nc{\BN}{{\mathbb{N}}}
\nc{\BP}{{\mathbb{P}}}
\nc{\BR}{{\mathbb{R}}}
\nc{\BZ}{{\mathbb{Z}}}

\nc{\Mod}{{{\mathcal M}od}}

\nc{\CA}{{\mathcal{A}}}
\nc{\CB}{{\mathcal{B}}}
\nc{\D}{{\mathcal{D}}}
\nc{\CE}{{\mathcal{E}}}
\nc{\CF}{{\mathcal{F}}}
\nc{\CG}{{\mathcal{G}}}
\nc{\CL}{{\mathcal{L}}}
\nc{\CM}{{\mathcal{M}}}
\nc{\CN}{{\mathcal{N}}}
\nc{\CO}{{\mathcal{O}}}
\nc{\CP}{{\mathcal{P}}}
\nc{\CQ}{{\mathcal{Q}}}
\nc{\CS}{{\mathcal{S}}}
\nc{\CT}{{\mathcal{T}}}
\nc{\CU}{{\mathcal{U}}}
\nc{\CV}{{\mathcal{V}}}
\nc{\CW}{{\mathcal{W}}}
\nc{\IC}{{\mathcal{IC}}}

\nc{\cM}{{\check{\mathcal M}}{}}
\nc{\csM}{{\check{\mathcal A}}{}}
\nc{\oM}{{\overset{\circ}{\mathcal M}}{}}
\nc{\obM}{{\overset{\circ}{\mathbf M}}{}}
\nc{\oCA}{{\overset{\circ}{\mathcal A}}{}}
\nc{\obA}{{\overset{\circ}{\mathbf A}}{}}
\nc{\ooM}{{\overset{\circ}{M}}{}}
\nc{\osM}{{\overset{\circ}{\mathsf M}}{}}
\nc{\vM}{{\overset{\bullet}{\mathcal M}}{}}
\nc{\nM}{{\underset{\bullet}{\mathcal M}}{}}
\nc{\oD}{{\overset{\circ}{\mathcal D}}{}}
\nc{\obD}{{\overset{\circ}{\mathbf D}}{}}
\nc{\oA}{{\overset{\circ}{\mathbb A}}{}}
\nc{\op}{{\overset{\bullet}{\mathbf p}}{}}
\nc{\cp}{{\overset{\circ}{\mathbf p}}{}}
\nc{\oU}{{\overset{\bullet}{\mathcal U}}{}}

\nc{\fa}{{\mathfrak{a}}}
\nc{\fb}{{\mathfrak{b}}}
\nc{\fg}{{\mathfrak{g}}}
\nc{\fgl}{{\mathfrak{gl}}}
\nc{\fj}{{\mathfrak{j}}}
\nc{\fn}{{\mathfrak{n}}}
\nc{\fu}{{\mathfrak{u}}}
\nc{\fp}{{\mathfrak{p}}}
\nc{\fr}{{\mathfrak{r}}}
\nc{\fsl}{{\mathfrak{sl}}}
\nc{\hsl}{{\widehat{\mathfrak{sl}}}}
\nc{\hgl}{{\widehat{\mathfrak{gl}}}}
\nc{\hg}{{\widehat{\mathfrak{g}}}}
\nc{\chg}{{\widehat{\mathfrak{g}}}{}^\vee}
\nc{\hn}{{\widehat{\mathfrak{n}}}}
\nc{\chn}{{\widehat{\mathfrak{n}}}{}^\vee}

\nc{\fA}{{\mathfrak{A}}}
\nc{\fB}{{\mathfrak{B}}}
\nc{\fD}{{\mathfrak{D}}}
\nc{\fE}{{\mathfrak{E}}}
\nc{\fF}{{\mathfrak{F}}}
\nc{\fG}{{\mathfrak{G}}}
\nc{\fK}{{\mathfrak{K}}}
\nc{\fL}{{\mathfrak{L}}}
\nc{\fM}{{\mathfrak{M}}}
\nc{\fN}{{\mathfrak{N}}}
\nc{\fP}{{\mathfrak{P}}}
\nc{\fU}{{\mathfrak{U}}}

\nc{\bb}{{\mathbf{b}}}
\nc{\bc}{{\mathbf{c}}}
\nc{\be}{{\mathbf{e}}}
\nc{\bj}{{\mathbf{j}}}
\nc{\bn}{{\mathbf{n}}}
\nc{\bp}{{\mathbf{p}}}
\nc{\bq}{{\mathbf{q}}}
\nc{\bu}{{\mathbf{u}}}
\nc{\bv}{{\mathbf{v}}}
\nc{\bx}{{\mathbf{x}}}
\nc{\by}{{\mathbf{y}}}
\nc{\bw}{{\mathbf{w}}}
\nc{\bA}{{\mathbf{A}}}
\nc{\bB}{{\mathbf{B}}}
\nc{\bC}{{\mathbf{C}}}
\nc{\bD}{{\mathbf{D}}}
\nc{\bH}{{\mathbf{H}}}
\nc{\bM}{{\mathbf{M}}}
\nc{\bV}{{\mathbf{V}}}
\nc{\bW}{{\mathbf{W}}}
\nc{\bX}{{\mathbf{X}}}

\nc{\sA}{{\mathsf{A}}}
\nc{\sB}{{\mathsf{B}}}
\nc{\sD}{{\mathsf{D}}}
\nc{\sF}{{\mathsf{F}}}
\nc{\sK}{{\mathsf{K}}}
\nc{\sM}{{\mathsf{M}}}
\nc{\sO}{{\mathsf{O}}}
\nc{\sQ}{{\mathsf{Q}}}
\nc{\sP}{{\mathsf{P}}}
\nc{\sfp}{{\mathsf{p}}}
\nc{\sr}{{\mathsf{r}}}

\nc{\BK}{{\bar{K}}}

\nc{\tA}{{\widetilde{\mathbf{A}}}}
\nc{\TG}{{\tilde{G}}}
\nc{\TM}{{\widetilde{\mathbb{M}}}{}}
\nc{\tO}{{\widetilde{\mathsf{O}}}{}}
\nc{\TZ}{{\tilde{Z}}}
\nc{\tx}{{\tilde{x}}}
\nc{\tbv}{{\tilde{\bv}}}
\nc{\tfP}{{\widetilde{\mathfrak{P}}}{}}
\nc{\tz}{{\tilde{\zeta}}}
\nc{\tmu}{{\tilde{\mu}}}

\nc{\urho}{\underline{\rho}}
\nc{\uB}{\underline{B}}
\nc{\uC}{{\underline{\mathbb{C}}}}
\nc{\ui}{\underline{i}}
\nc{\uj}{\underline{j}}
\nc{\ofP}{{\overline{\mathfrak{P}}}}

\nc{\eps}{\varepsilon}
\nc{\hrho}{{\hat{\rho}}}

\nc{\one}{{\mathbf{1}}}
\nc{\two}{{\mathbf{t}}}

\nc{\Rep}{{\mathop{\operatorname{\rm Rep}}}}
\nc{\Sym}{{\mathop{\operatorname{\rm Sym}}}}
\nc{\Tot}{{\mathop{\operatorname{\rm Tot}}}}
\nc{\Spec}{{\mathop{\operatorname{\rm Spec}}}}
\nc{\Ker}{{\mathop{\operatorname{\rm Ker}}}}
\nc{\Hilb}{{\mathop{\operatorname{\rm Hilb}}}}
\nc{\End}{{\mathop{\operatorname{\rm End}}}}
\nc{\Ext}{{\mathop{\operatorname{\rm Ext}}}}
\nc{\Hom}{{\mathop{\operatorname{\rm Hom}}}}
\nc{\CHom}{{\mathop{\operatorname{{\mathcal{H}}\it om}}}}
\nc{\GL}{{\mathop{\operatorname{\rm GL}}}}
\nc{\gr}{{\mathop{\operatorname{\rm gr}}}}
\nc{\Id}{{\mathop{\operatorname{\rm Id}}}}
\nc{\rk}{{\mathop{\operatorname{\rm r}}}}
\nc{\de}{{\mathop{\operatorname{\rm def}}}}
\nc{\length}{{\mathop{\operatorname{\rm length}}}}
\nc{\supp}{{\mathop{\operatorname{\rm supp}}}}

\nc{\Bun}{{\mathsf{Bun}}}
\nc{\Cliff}{{\mathsf{Cliff}}}
\nc{\Gr}{{\mathsf{Gr}}}
\nc{\Fl}{{\mathsf{Fl}}}
\nc{\Fib}{{\mathsf{Fib}}}
\nc{\Coh}{{\mathsf{Coh}}}
\nc{\FCoh}{{\mathsf{FCoh}}}

\nc{\reg}{{\text{\rm reg}}}

\nc{\cplus}{{\mathbf{C}_+}}
\nc{\cminus}{{\mathbf{C}_-}}
\nc{\cthree}{{\mathbf{C}_*}}
\nc{\Qbar}{{\bar{Q}}}

\nc{\bh}{{\bar{h}}}
\nc{\bOmega}{{\overline{\Omega}}}

\nc{\seq}[1]{\stackrel{#1}{\sim}}

\title{Uhlenbeck spaces for $\BA^2$ and affine Lie algebra $\hsl_n$}

\author{Michael Finkelberg}
\address{Independent Moscow University,
Bolshoj Vlasjevskij Pereulok, dom 11,
Moscow 121002 Russia}
\email{fnklberg@@mccme.ru}

\author{Dennis Gaitsgory}
\address{Department of Mathematics,
The University of Chicago, Chicago,
IL 60637, USA}
\email{gaitsgde@@math.uchicago.edu}

\author{Alexander Kuznetsov}
\address{
Institute for Problems of Information Transmission,
 Russian Academy of Sciences,
19 Bolshoi Karetnyi, Moscow 101447, Russia
}
\email{sasha@@kuznetsov.mccme.ru}

\begin{document}

\begin{abstract}
We introduce an Uhlenbeck closure of the space of based maps from projective
line to the Kashiwara flag scheme of an untwisted affine Lie algebra.
For the algebra $\widehat{sl}_n$ this space of based maps is isomorphic to
the moduli space of locally free parabolic sheaves on $P^1\times P^1$
trivialized at infinity. The Uhlenbeck closure admits a resolution of
singularities: the moduli space of torsion free parabolic sheaves
on $P^1\times P^1$ trivialized at infinity. We compute the Intersection
Cohomology sheaf of the Uhlenbeck space using this resolution of singularities.
The moduli spaces of parabolic sheaves of various degrees are connected by
certain Hecke correspondences. We prove that these correspondences define
an action of $\widehat{sl}_n$ in the cohomology of the above moduli spaces.
\end{abstract}

\maketitle

\section{Introduction}

\subsection{} For a symmetrizable Cartan matrix $A$, and the corresponding
Kac-Moody algebra $\fg(A)$, M.~Kashiwara has introduced a remarkable {\em flag
scheme} $\CB(A)$ ~\cite{k}. It shares many properties of the usual flag
varieties of semisimple Lie algebras. For one thing, if $\bC$ is a smooth
projective curve of genus 0, and $\bc\in\bC$ a marked point, the space
$\obM^\alpha(A)$ of {\em based maps} from $(\bC,\bc)$ to $(\CB(A),B_0)$ of
degree $\alpha$ turns out surprisingly to be a smooth finite-dimensional
quasiaffine variety, though $\CB(A)$ itself is of infinite type in general.

In case $\fg(A)$ is semisimple, V.~Drinfeld has introduced a remarkable
affine closure $\bM^\alpha(A)\supset\obM^\alpha(A)$ (the space of {\em based
quasimaps}, alias {\em Zastava} space) which has found applications in the
study of quantum groups at roots of unity and geometric Eisenstein series.
In fact, Drinfeld's definition works for arbitrary symmetrizable $A$, but
$\bM^\alpha(A)$ turns out to be of infinite type in general.

However, if $\fg(A)$ is an {\em untwisted affine} Lie algebra, it appears that
$\bM^\alpha(A)$ possesses a partial resolution $\fM^\alpha(A)\to\bM^\alpha(A)$
with rather favorable properties (for one thing, $\fM^\alpha(A)$ is of finite
type). The construction of $\fM^\alpha(A)$ we propose at the moment is
quite cumbersome, in a sense it occupies the bulk of this paper. As a shortcut
we conjecture that $\fM^\alpha(A)$ is {\em normal}; it would imply that it is
just the normalization of $\bM^\alpha(A)$.

\subsection{} In case $A=\tA_{n-1}$ we have $\fg(A)=\hsl_n$, and
$\fM^\alpha:=\fM^\alpha(\tA_{n-1})$ admits a semismall resolution of
singularities $\varpi_\alpha:\ \CM^\alpha\to\fM^\alpha$. The quasiprojective
variety $\CM^\alpha$ is not new; it is just the moduli space of
torsion free {\em parabolic sheaves} of degree $\alpha$
on the surface $\bC\times\BP^1$ trivialized at infinity. For a torsion free
parabolic sheaf $\CF_\bullet\in\CM^\alpha$ one can define its {\em saturation}
$\CN(\CF)_\bullet$ which is a {\em locally free} parabolic sheaf containing
$\CF_\bullet$, and {\em defect} $\de(\CF_\bullet)$ which is roughly speaking
a colored zero-cycle on $\bC\times\BP^1$ measuring the quotient
$\CN(\CF)_\bullet/\CF_\bullet$. The proper map
$\varpi_\alpha:\ \CM^\alpha\to\fM^\alpha$ glues together various
parabolic sheaves with the same saturation and defect.

Thus, the closure $\fM^\alpha\supset\obM^\alpha$ may be obtained in two steps.
First we view $\obM^\alpha$ as the moduli space of
locally free parabolic sheaves of degree $\alpha$
on the surface $\bC\times\BP^1$ trivialized at infinity, and put it inside
the moduli space of torsion free parabolic sheaves $\CM^\alpha$. Second,
we glue together certain torsion free parabolic sheaves. This idea is not
new: for the moduli spaces of vector bundles on surfaces it gives rise to
{\em Uhlenbeck} compactifications. This is why we call $\fM^\alpha$ an
{\em Uhlenbeck flag space} for $\BA^2$ (though we work with the surface
$\bC\times\BP^1$, the trivialization at infinity essentially leaves us with
$\BA^2\subset\bC\times\BP^1$).

Unfortunately, the ``Uhlenbeck compactification'' of a moduli space of
vector bundles on a surface has been given a rigorous algebraic geometric
meaning
only in few instances, notably for the vector bundles on $\bC\times\BP^1$
trivialized at infinity (or equivalently, 
vector bundles on $\BP^2$ trivialized
at infinity) in the remarkable works of H.~Nakajima on Quiver varieties.

So roughly speaking we cook up our $\fM^\alpha$ from Drinfeld's Zastava
space $\bM^\alpha$ and Nakajima's Uhlenbeck space for $\BA^2$. Though
the exposition in the main body of the paper concerns the case
$\fg(A)=\hsl_n$, we spell it out in such a way that the construction
carries out without changes for an arbitrary untwisted affine Lie algebra
$\fg(A)$, cf. ~\ref{speculation}.

\subsection{}
The particularity of the case $\fg(A)=\hsl_n$ lies in the existence of
a semismall resolution of singularities
$\varpi_\alpha:\ \CM^\alpha\to\fM^\alpha$. This is similar to the existence
of a {\em small resolution}, due to G.~Laumon, of the Zastava space
$\bM^\alpha(\fsl_n)$. We apply the resolution $\varpi_\alpha$ to compute
the Intersection Cohomology sheaf $\IC(\fM^\alpha)$, similarly to ~\cite{ku},
where Laumon's resolution was used to compute $\IC(\bM^\alpha(\fsl_n))$.
The necessary information about the fibers of $\varpi_\alpha$ was already
obtained in ~\cite{fk2}, ~\cite{nk}, so in a sense, all the hard work was
already done a long time ago. The generating function of the $\IC$-stalks
is governed by the product of Kostant partition function for $\hsl_n$,
and another partition function, arising from the invariants of a principal
nilpotent element of $\fsl_n$ in the nilpotent radical of the maximal
parabolic subalgebra of $\hsl_n$.

For an arbitrary untwisted affine $\fg(A)$ we propose a conjectural answer
for the stalks of $\IC(\fM^\alpha(A))$ in ~\ref{conj}.

\subsection{}
We also study another moduli space $\CM^\alpha_{\mathbf{gt}}\supset\CM^\alpha$
of parabolic torsion free sheaves of degree $\alpha$ on $\bC\times\BP^1$,
where we relax the condition of triviality at infinity, and impose only a
condition that a torsion free sheaf $\CF_0$ is {\em generically trivial},
that is trivial on some line $c\times\BP^1$.
For any $\alpha,\gamma$ there is a closed subvariety of middle
dimension ({\em Hecke correspondence}) $\fE_\alpha^\gamma\subset
\CM_{\mathbf{gt}}^\alpha\times\CM_{\mathbf{gt}}^{\alpha+\gamma}$.
It is formed by pairs of parabolic sheaves 
such that the second one is a subsheaf
of the first one. The top-dimensional irreducible components of
$\fE_\alpha^\gamma$ are naturally numbered by the isomorphism classes
$\kappa\in\fK(\gamma)$ of $\gamma$-dimensional nilpotent representations
of the cyclic quiver $\tA_{n-1}$, independently of $\alpha$. For
$\kappa\in\fK(\gamma)$ the corresponding irreducible component 
$\fE_\alpha^\kappa$,
viewed as a correspondence between $\CM_{\mathbf{gt}}^\alpha$ and 
$\CM_{\mathbf{gt}}^{\alpha+\gamma}$,
defines two operators:
$$e_\kappa:\ H^\bullet(\CM_{\mathbf{gt}}^\alpha)\rightleftharpoons
H^\bullet(\CM_{\mathbf{gt}}^{\alpha+\gamma})\ :f_\kappa$$
Let $\bH$ denote the generic Hall algebra of nilpotent representations of
the cyclic quiver $\tA_{n-1}$ at $q=1$. It turns
out that the linear span of operators $e_\kappa$ is closed under composition;
the algebra they form is naturally isomorphic to $\bH$, and the isomorphism
takes $e_\kappa$ to the element of $\bH$ corresponding to the isomorphism
class $\kappa$. Moreover, for the isomorphism classes of simple representations
$\kappa=\{\{i\}\},\ i\in\BZ/n\BZ$, the corresponding operators $e_i,f_i$
define the action of the Chevalley generators of $\fg(\tA_{n-1})=\hsl_n$
on $\bigoplus_\alpha H^\bullet(\CM_{\mathbf{gt}}^\alpha)$. This $\hsl_n$-action
has central charge 2. This is a partial realization of the programme outlined
in ~\cite{fk2} ~1.3.

\subsection{} Let us say a few words about the structure of the paper.
In \S2 we recall the well known facts about
the Kashiwara flag scheme for $\hsl_n$, and various realizations thereof.
In \S3 we introduce the moduli space of torsion free parabolic sheaves
$\CM^\alpha$, and construct a family of regular functions on it, which
will be used in the definition of the resolution
$\varpi_\alpha:\ \CM^\alpha\to\fM^\alpha$.
In \S4 we recall the Drinfeld's spaces of based maps and quasimaps
$\obM^\alpha\subset\bM^\alpha$, and define the Uhlenbeck flag space
$\fM^\alpha$ as a closure of $\obM^\alpha$ in some quasiaffine embedding
(like Schubert varieties are closures of Schubert cells in the usual flag
varieties).
In \S5 we construct the resolution $\varpi_\alpha$, and in \S6 we
compute $\IC(\fM^\alpha)$. Note that while the generating function of
the stalks of $\varpi_{\alpha*}\IC(\CM^\alpha)$ involves the Kostant
partition function of $\hgl_n$, the generating function of the $\IC$-stalks
of $\fM^\alpha$ involves the Kostant partition function of $\hsl_n$:
the {\em semi}smallness of $\varpi_\alpha$ kills the extra imaginary roots.
In \S7 we study the Hecke correspondences; among other things, they are used
in the proof of connectedness of $\CM^\alpha$.

\subsection{} Our main motivation was to understand the algebraic geometric
meaning of Uhlenbeck compactifications. We did not really succeed (for one
thing, we are bound to the surface $\BA^2$ with fixed coordinates); the
present work may be viewed just as an indication what to look for.
We benefited strongly from the explanations by V.~Baranovsky, V.~Drinfeld
and V.~Ginzburg about Uhlenbeck compactifications. Moreover, this work owes
its very existence to the ideas and suggestions of V.Drinfeld. We are also
grateful to O.Schiffmann for bringing the reference ~\cite{fm} to our
attention. In the course of our
study of Uhlenbeck spaces, M.F. has enjoyed the hospitality and support
of the IHES, the Universit\'e Cergy-Pontoise, the Hebrew University of
Jerusalem, and the University of Chicago. His research was conducted for the
Clay Mathematical Institute. D.G. is a Prize Fellow of the Clay Mathematical
Institute. A.K. was partially supported by RFFI grants
99-01-01144 and 99-01-01204.

\section{Kashiwara flag scheme for $\hsl_n$}

\subsection{} Recall that the affine Lie algebra $\hsl_n$
is the canonical central extension
$$0\to\BC\to\hsl_n\to\fsl_n\otimes\BC((t^{-1}))\to0$$
Let us fix an $n$-dimensional vector space $V$ with a basis $v_1,\ldots,v_n$,
and identify $\fsl_n$ with $\fsl(V)$. For $1\leq i\ne j\leq n$ we denote by
$E_{ij}\in\fsl_n$ the operator taking $v_j$ to $v_i$, and annihilating other
base vectors. Then the Chevalley generators of $\hsl_n$ are as follows:
$e_0=t^{-1}E_{n1},\ f_0=tE_{1n},\ h_0=[e_0,f_0]$; for $1\leq i\leq n-1$ we set
$e_i=E_{i,i+1},\ f_i=E_{i+1,i},\ h_i=[e_i,f_i]$. Thus the simple positive
coroots are naturally numbered by $0\leq i\leq n-1$. We will identify this set
with $I:=\BZ/n\BZ$. We will denote by $Y$ the coroot lattice $\BZ[I]$,
and by $X$ the dual weight lattice. We denote the perfect pairing
$X\times Y\to\BZ$ by $\langle,\rangle$,
and the basis of $X$ dual to $I$ consists
of fundamental weights $\omega_i,\ i\in I$. Thus $\langle\omega_j,i\rangle=
\delta_{ij}$. A simple root dual to a simple coroot $i\in I$ will be denoted
by $i'\in X$. 

For a dominant weight $X^+\ni\lambda=\sum_Il_i\omega_i,\ l_i\in\BN$, we
denote by $V_\lambda$ the corresponding highest weight integrable
$\hsl_n$-module (its highest vector is annihilated by $e_i,\ i\in I$).
We denote by $V_\lambda^*$ the dual (pro-finite dimensional) vector space.

\subsection{Fundamental representations}
\label{wedge}
The general reference for this subsection is ~\cite{bbe}.
Recall the semi-infinite wedge construction of the fundamental representations
$V_{\omega_i}$. Let $\bV$ denote the Tate vector space
$V\otimes\BC((t^{-1}))$. Then the Tate vector space $\bW:=\bV\oplus\bV^*$
has a natural symmetric
bilinear form which gives rise to the Clifford algebra $\Cliff(\bW)$.
We choose a compact lattice $L_V=V\otimes\BC[[t^{-1}]]\subset\bV$, and consider
a compact lattice $L_W=L_V\oplus L_V^\perp\subset\bW$. Then $L_W$ is an
isotropic subspace of $\bW$, and its exterior algebra embeds naturally
into the Clifford algebra: $\Lambda^\bullet(L_W)\subset\Cliff(\bW)$.

We define the Clifford module $\CQ$ as
$\operatorname{Ind}_{\Lambda^\bullet(L_W)}^{\Cliff(\bW)}$. In fact, its
isomorphism class is independent of the choice of compact lattice
$L_V\subset\bV$.

Consider an arbitrary compact lattice $L_1\subset\bV$ and another compact
lattice $L_2\subset L_1^\perp\subset\bV^*$. We set
$L_{1,2}:=L_1\oplus L_2\subset\bW$. Then the invariants $\CQ^{L_{1,2}}$
form a finite
dimensional vector subspace canonically
isomorphic to $\Lambda^*(L_2^\perp/L_1)\otimes\det(L_2^\perp)$.
Here $\Lambda^*(?)$ is a vector space dual to the exterior algebra
$\Lambda^\bullet(?)$, and $\det(L_2^\perp)$ is the relative determinant
of the lattice $L_2^\perp\subset\bV$ with respect to $L_V$.
Clearly, $\CQ$ is a union of $\CQ^{L_{1,2}}$ as $L_1, L_2$ shrink.

The algebra $\hgl_n\supset\hsl_n$ acts naturally on $\CQ$.
It is well known that for any $i\in I$ there is a canonical embedding
$s_i:\ V_{\omega_i}\to\CQ$. In fact, $\CQ$ is the direct sum of fundamental
representations of $\hgl_n$.

\subsection{Pl\"ucker equations}
\label{Pluck}
Kashiwara ~\cite{k} defines the flag scheme $\CB$ for $\hsl_n$ as the
(infinite type) subscheme of $\prod_{\lambda\in X^+}\BP(V_\lambda^*)$
cut out by {\em Pl\"ucker equations}:

A collection of lines $(\ell_\lambda\subset V_\lambda^*)_{\lambda\in X^+}$
satisfies Pl\"ucker equations if

(a) For any nonzero $\hsl_n$-morphism $\varphi:\ V_\lambda^*\hat{\otimes}
V_\mu^*\to V_{\lambda+\mu}^*$ we have $\varphi(\ell_\lambda\otimes\ell_\mu)=
\ell_{\lambda+\mu}$;

(b) For any $\hsl_n$-morphism $\varphi:\ V_\lambda^*\hat{\otimes}
V_\mu^*\to V_\nu^*$ such that $\nu<\lambda+\mu$
we have $\varphi(\ell_\lambda\otimes\ell_\mu)=0$.

\medskip

The inverse image of the line bundle $\CO(1)$ on $\BP(V_\lambda^*)$ is
the line bundle on $\CB$ denoted by $\CL_\lambda$.
We have $\Gamma(\CB,\CL_\lambda)=V_\lambda$.

Note that the Pl\"ucker equation (a) above implies that $\CB$ embeds
as a closed subscheme into $\prod_{i\in I}\BP(V_{\omega_i}^*)$.

\subsection{Discrete lattices}
\label{flag}
The above definition works in the generality of an arbitrary symmetrizable
Kac-Moody algebra. In the particular case of $\hsl_n$ there is another
well known definition of $\CB$ in terms of periodic flags in the Tate
vector space $\bV$. Namely, $\CB$ is a scheme (of infinite type) parametrizing
collections of discrete lattices $(F_k\subset\bV)_{k\in\BZ}$ such that

(a) The kernel and cokernel of the natural map
$F_0\oplus V\otimes\BC[[t^{-1}]]\to\bV$ have the same dimension.

(b) $F_k\subset F_{k+1}$, and $\dim(F_{k+1}/F_k)=1$ for any $k$;

(c) $F_{k+n}=t^{-1}F_k$ for any $k$;

\medskip

Let us construct an isomorphism from the second definition of $\CB$ to
the first one. To this end let us temporarily denote $\CB$ in the first
(resp. second) definition by $\CB_1$ (resp. $\CB_2$).
Given a flag $(F_k)$ and $0\leq i\leq n-1$, we consider
a discrete lattice $F_{W,i}=F_i\oplus F_i^\perp\subset\bW$. It is well known
that the coinvariants of the Clifford module $\CQ_{F_{W,i}}$ are
one-dimensional. Composing with the canonical embedding
$s_i:\ V_{\omega_i}\to\CQ$ (see ~\ref{wedge}) we obtain a projection
$p_i:\ V_{\omega_i}\to\CQ_{F_{W,i}}$.

Thus we have constructed a line bundle $\CL_i$ over $\CB_2$
(with a fiber over $(F_k)_{k\in\BZ}$ equal to $\CQ_{F_{W,i}}$) together
with a surjection $p_i:\ V_{\omega_i}\otimes\CO_{\CB_2}\twoheadrightarrow
\CL_i$. It defines a map $\CB_2\to\BP(V_{\omega_i}^*)$, and taking the
product over $i\in I$ we obtain an embedding
$\CB_2\hookrightarrow\prod_{i\in I}\BP(V_{\omega_i}^*)$ which identifies
it with the image of Pl\"ucker embedding of $\CB_1$.
This way $\CL_i$ on $\CB_2$ gets identified with $\CL_{\omega_i}$ on $\CB_1$.

\subsection{Parabolic vector bundles}
\label{P1 moduli}
The second definition of the flag scheme $\CB$ translates immediately into
the language of vector bundles on $\BP^1$. Namely, let $\bX$ be a smooth
projective curve of genus 0. We choose two distinct points $\by,\bx\in \bX$
and a global rational coordinate $t:\ \bX\to\BP^1$ such that $t(\by)=0,\
t(\bx)=\infty$.

Then $\CB_2$ is isomorphic to the moduli space $\CB_3$ of {\em parabolic
vector bundles} on $\bX$ with a trivialization in the formal neighbourhood
of $\bx\in \bX$. More precisely, we consider the moduli space of the collections
$(\fF_k,\tau)_{k\in\BZ}$ where

(a) $\fF_k$ is a vector bundle on $\bX$ of degree $k$ and rank $n$;

(b) $\fF_k\subset\fF_{k+1}$, and $\fF_{k+1}/\fF_k$ 
is supported at $\by\in \bX$
for any $k$;

(c) $\fF_{k+n}=\fF_k(\by)$ for any $k$;

(d) $\tau$ is a trivialization of $\fF_0$ restricted to the formal
neighbourhood $\bX_{\widehat \bx}$ of $\bx\in \bX$ (and hence $\tau$ is a
trivialization of any $\fF_k$ in $\bX_{\widehat \bx}$).

\medskip

Let us recall the isomorphism from $\CB_3$ to $\CB_2$.
Given $(\fF_k,\tau)_{k\in\BZ}$ we define the flag of discrete lattices
$(F_k)_{k\in\BZ}$ as follows. Our coordinate $t:\ \bX\to\BP^1$ identifies
$\CO_{\bX_{\widehat \bx}}$ with $\BC[[t^{-1}]]$. Hence $\tau$ identifies
$\fF_k|_{\bX_{\widehat \bx}}$
with $V\otimes\BC[[t^{-1}]]$. Now the space of global
sections $\Gamma(\bX-\bx,\fF_k)$ embeds as a discrete lattice $F_k$ into
$\Gamma(\bX_{\widehat \bx}-\bx,\fF_k)=V\otimes\BC((t^{-1}))$. One checks easily
that the conditions ~(a--c) above imply the conditions ~\ref{flag} ~(a--c).

Under this isomorphism, the fiber of the line bundle $\CL_i$ at a point
$(F_k)_{k\in\BZ}$ gets identified with the determinant of cohomology
$\det R\Gamma(\bX,\fF_i)$.

\subsection{Schubert divisors}
\label{Delta}
Recall that $V$ is a vector space with a basis $v_1,\ldots,v_n$. We define
a complete flag of vector subspaces $0=V_0\subset V_1\subset\ldots\subset
V_{n-1}\subset V_n=V$ where $V_i=\langle v_1,\ldots,v_i\rangle$.
We define a transversal flag $V=V^0\supset V^{-1}\supset\ldots\supset
V^{1-n}\supset V^{-n}=0$ where $V^{-j}=\langle v_n,\ldots,v_{j+1}\rangle$.
We denote by $\sB$ the flag variety of $\fsl_n$. So we have two distinguished
points $V^\bullet, V_\bullet\in\sB$. Let
$\fb\subset\fsl_n$ be a Borel subalgebra formed by all the operators preserving
our flag $V_\bullet$. Let $\bB\subset SL_n$ be the corresponding Borel
subgroup. Let $\fn\subset\fb$ be the nilpotent radical.

We define a subalgebra $\hat\fb\subset\hsl_n$ as a full preimage in the
central extension of a subalgebra $\fb\oplus\fsl_n\otimes t^{-1}\BC[[t^{-1}]]
\subset\fsl_n\otimes\BC((t^{-1}))$. Let $\hat\bB$ be the corresponding
proalgebraic group. According to ~\cite{k}, $\hat\bB$ acts on $\CB$ with
a unique open orbit $\CU\subset\CB$. The complement $\CB-\CU$ is a union of
$n+1$ irreducible Cartier divisors naturally numbered by $I:\ \CB-\CU=
\bigcup_{i\in I}\Delta_i$. We have $\CL_i=\CL_{\omega_i}=\CO(\Delta_i)$,
see ~\cite{k2}.

Finally, recall that $\Delta_0$ is cut out by
the condition that $\CF_0$ is a nontrivial vector
bundle on $\bX$.

\subsection{Base point}
\label{base point}
We choose a base point $B_0\in\CU\subset\CB$ as follows. In the setup
of ~\ref{Pluck} we set $B_0=(\ell^0_\lambda)_{\lambda\in X^+}$ where
$\ell^0_\lambda$ is the unique line in $V_\lambda^*$ killed by all
$f_i,\ i\in I$. Equivalently, in the setup of ~\ref{flag} we have
$B_0=(F_k)_{k\in\BZ}$ where for $-n\leq k\leq 0$ we have
$F_k=V\otimes t\BC[t]\oplus V^k$. Equivalently, in the setup
of ~\ref{P1 moduli} we have $B_0=(\fF_k,\tau)_{k\in\BZ}$ where
$\fF_0=V\otimes\CO_\bX,\ \tau$ is the tautological trivialization, and
for $-n\leq k\leq0$ the local sections of $\fF_k$ are those sections of
$\fF_0=V\otimes\CO_\bX$ which take value in $V^k$ at $\by\in \bX$.

\subsection{Beilinson-Drinfeld-Kottwitz flags}
\label{BDK}
We recall the construction ~\cite{g} of an ind-scheme of ind-finite type
``approximating'' the infinite type scheme $\CB$. For a positive integer
$a$ let $\bX^{(a)}$ denote the $a$-th symmetric power of $\bX$. For a test
scheme $S$, and an $S$-point $y$ of $\bX^{(a)}$, we may view the graph
$\Gamma_y$ of $y$ as a subscheme of $S\times\bX$ (finite over $S$).

Following ~\cite{g}, we define the ind-scheme $\fB^a$ representing the
functor associating to a test scheme $S$ the set of quadruples
$(y,\CV,\varsigma,\CV^\bullet_\by)$ where

$y$ is an $S$-point of $(\bX-\bx)^{(a)}$;

$\CV$ is an $SL_n$-bundle on $S\times\bX$;

$\varsigma$ is a trivialization $\CV|_{S\times\bX-\Gamma_y}\to
V\otimes\CO_{S\times\bX-\Gamma_y}$;

$\CV^\bullet_\by$ is a reduction of $\CV|_{S\times\by}$ to $\bB\subset SL_n$.

\medskip

$\fB^a$ is equipped with an evident projection $p_a:\ \fB^a\to(\bX-\bx)^{(a)}$,
and with a section $s_a:\ (\bX-\bx)^{(a)}\to\fB^a$ defined as follows.
For $s_a(y)\in\fB^a$ we have: $\CV=V\otimes\CO_{S\times\bX}$ is a trivial
$SL_n$-bundle; $\varsigma=\Id$ is the tautological trivialization;
$\CV^\bullet_\by$ is given by a constant flag
$V^\bullet\otimes\CO_{S\times\by}$ in
$\CV_{S\times\by}=V\otimes\CO_{S\times\by}$.

We have an evident morphism $m_a:\ \fB^a\to\CB$ restricting a rational
trivialization $\varsigma$ to the formal neighbourhood $\bX_{\widehat\bx}$
of $\bx$ in $\bX$. Note that $m_a$ contracts the section $s_a((\bX-\bx)^{(a)})$
to the base point $B_0\in\CB$.

\subsection{Kashiwara Grassmannian}
\label{kg} 
Kashiwara scheme $\CB$ has an important parabolic version $\CG$ which we
presently recall. In the setup of ~\ref{Pluck}, for $i=0$, the line bundle
$\CL_{\omega_0}$ on $\CB$ defines a morphism from $\CB$ to 
$\BP(V_{\omega_0}^*)$, and $\CG$ is the image of this morphism. We have a
fiber bundle $\CB\to\CG$ with the typical fiber $\sB$. Thus, the line
bundle $\CL_{\omega_0}$ on $\CB$ descends to the ample 
{\em determinant line bundle} $\CL_0$ on $\CG$. 

Equivalently, in the setup of ~\ref{flag}, $\CG$ is the moduli scheme
of discrete lattices $F\subset\bV$ satisfying
the condition (a) of {\em loc. cit.}, such that $F\subset t^{-1}F$.

Equivalently, in the setup of ~\ref{P1 moduli}, $\CG$ is the moduli scheme
of pairs $(\fF,\tau)$ where $\fF$ is an $SL_n$-bundle on $\bX$, and
$\tau$ is a trivialization of $\fF$ in the formal neighbourhood of 
$\bx\in\bX$.

We have a divisor $\Delta_0\subset\CG$ cut out by the condition that 
$\CF_0$ is a nontrivial vector bundle on $\bX$, and $\CL_0=\CO(\Delta_0)$.
Also, we have a base point $G_0\in\CG$ which is the image of $B_0\in\CB$.
Finally, in the setup of ~\ref{BDK}, 
for $a\in\BN$ we have the ind-scheme $\fG^a$ ({\em Beilinson-Drinfeld
Grassmannian}) representing the functor associating to a test scheme
$S$ the set of triples $(y,\CV,\varsigma)$ as in {\em loc. cit.}
We have an evident morphism $m_a:\ \fG^a\to\CG$.

\section{Parabolic sheaves on $\BA^2$}

\subsection{}
\label{surfaces}
Let $\bC$ be a smooth projective curve
of genus 0. We choose two distinct points $\bb,\bc\in \bC$
and a global rational coordinate $z:\ \bC\to\BP^1$ such that $z(\bb)=0,\
z(\bc)=\infty$.

We consider a smooth projective surface $\CS':=\bC\times \bX$ with a normal
crossing divisor $\bD':=\bC\times \bx\bigcup \bc\times \bX$. 
Note that $\CS'-\bD'$
is the affine plane $\BA^2$ with coordinates $z,t$.

Blowing up the point $\bc\times \bx\in\CS'$ we obtain a surface $\CS$ with
a projection $p:\ \CS\to\CS'$. It is well known that one can blow down
the proper transform of $\bD'$ in $\CS$ to obtain $q:\ \CS\to\CS''$.
The surfaces $\CS',\CS''$ have a common open subscheme
$\CS'\supset\BA^2\subset\CS''$, and the complement $\CS''-\BA^2$ is a
smooth divisor $\bD''\subset\CS''$. In fact, $\CS''$ is isomorphic to
$\BP^2$, and $\bD''$ is a projective line.

Finally, we introduce a divisor $\bD_0:=\bC\times \by\subset\CS'$. Note that
$\bD_0\cap\BA^2$ is cut out by the equation $t=0$.

\subsection{Torsion free sheaves}
\label{grassmann}
For a positive integer $a$ let $\CA^a\supset\oCA^a$ 
denote the fine moduli space of
torsion free (resp. locally free) 
coherent sheaves $\CF$ on $\CS'$ of rank $n$, and second
Chern class $a$, equipped with a trivialization at $\bD':\ \CF|_{\bD'}=
V\otimes\CO_{\bD'}$. Its existence is proved in ~\cite{hl}, and its
smoothness is well known. For the reader's convenience let us recall the
argument.

\begin{lemma}
\label{kuz}
$\CA^a$ is smooth.
\end{lemma}

{\sl Proof:} Let $\CF\in\CA^a$ be a torsion free sheaf. 
The obstruction to smoothness of $\CA^a$ at $\CF$ lies in 
$\Ext^2(\CF,\CF(-\bD'))$ which by Serre duality is a vector space dual to  
$\Hom(\CF,\CF(\bD')\otimes\Omega^2_{\CS'})\simeq\Hom(\CF,\CF(-\bD'))$.
We claim that the latter vector space is zero, which at the same time
proves that $\CF$ has no infinitesimal automorphisms.
%that is the moduli stack $\CA^a$ is in fact the fine moduli space. 
In effect, since $\CF|_{\bC\times\bx}$
is a trivial vector bundle on $\bC=\BP^1$, for a general $x\in\bX$ the
restriction $\CF_x:=\CF|_{\bC\times x}$ 
is also a trivial vector bundle on $\bC$.
But then $\Hom(\CF_x,\CF_x(-\bc))=0$ for a general $x\in\bX$, and hence
already $\Hom(\CF,\CF(-\bc\times\bX))=0$.
\qed

\medskip

We will use an equivalent definition of $\CA^a$ going back to ~\cite{a}.
Namely, let $\CA^a_1$ denote a fine moduli space of torsion free coherent
sheaves $\CE$ on $\CS''$ of rank $n$, and second Chern class $a$, equipped
with a trivialization at $\bD'':\  \CE|_{\bD''}=
V\otimes\CO_{\bD''}$. Its existence is proved in ~\cite{hl}.

Following ~\cite{a}, we construct an isomorphism $\xi_a$
from $\CA^a$ to $\CA^a_1$
sending $\CF$ to $\CE:=q_*p^*\CF$ (notations of ~\ref{surfaces});
the inverse isomorphism from $\CA^a_1$
to $\CA^a$ sends $\CE$ to $\CF:=p_*q^*\CE$.

\medskip

Given a torsion free sheaf $\CF\in\CA^a$ and a point $s\in\BA^2\subset\CS'$
we define a {\em saturation} at
$s:\ \CN_s(\CF):=\bj_{s*}\bj_s^*\CF\supset\CF$
where $\bj_s:\ \CS'-s\hookrightarrow\CS'$ is an open embedding. It is well
known that $\CN_s(\CF)$ is a torsion free sheaf locally free at $s$.
We define a {\em defect} at $s:\ \de_s(\CF)$ as the length of the torsion
sheaf $\CN_s(\CF)/\CF$. Finally, $\CN(\CF)$ denotes the total
{\em saturation} of $\CF$, that is, $\CN(\CF):=\bj_{S*}\bj_S^*\CF$
where $S\subset\BA^2\subset\CS'$
is a finite subset such that $\CF$ is locally free
off $S$, and $\bj_S:\ \CS'-S\hookrightarrow\CS'$ is an open embedding.
The sum $\sum_{s\in S}\de_s(\CF)\cdot s\in\Sym^d(\BA^2)$ is the total
{\em defect} $\de(\CF)$.

\subsection{Quiver description}
\label{nakajima}
Nakajima ~(\cite{n1}, Theorem 2.1) gives another equivalent definition
of $\CA(n,a)=\CA^a_1$ as a certain quiver variety.
Recall that $\CA(n,a)$ is a moduli space of certain
linear algebra data $(B_1,B_2,\imath,\jmath)$, see ~{\em loc. cit.}
Here $B_1,B_2\in\End(W)$ where $W=\BC^a,\ \imath\in\Hom(V,W),\
\jmath\in\Hom(W,V)$ satisfy a condition $[B_1,B_2]+\imath\jmath=0$.
Nakajima defines $\CF$ as the middle cohomology of
a certain monad on $\CS''$ constructed from these linear algebra data.

Recall that
$z$ is our coordinate on $\bC$ which identifies $\bC-\bc$ with $\BA^1$.
Also, $t$ is our coordinate on $\bX$ which identifies $\bX-\bx$ with $\BA^1$.
The restriction of this monad to $\BA^2\subset\CS''$ looks as follows:
\begin{equation}
\label{monad}
0\to W\otimes\CO_{\BA^2}\stackrel{\fa}{\longrightarrow}
(W\oplus W\oplus V)\otimes\CO_{\BA^2}\stackrel{\fb}{\longrightarrow}
W\otimes\CO_{\BA^2}\to0
\end{equation}
where $\fa$ sends $w\in W\otimes\CO_{\BA^2}$ to
$((B_1-z)w,(B_2-t)w,\jmath w)$, and $\fb$ sends a triple $(w_1,w_2,v)$ to
the sum $-(B_2-t)w_1+(B_1-z)w_2+\imath v$.

\subsection{Parabolic sheaves}
\label{parab}
Let $\alpha=\sum_{i\in I}a_ii\in\BN[I]\subset Y$ be a positive coroot
combination. A {\em parabolic sheaf} $\CF_\bullet$ of degree $\alpha$
on $\BA^2$ is an infinite flag of torsion free coherent sheaves of rank $n$ on
$\CS':\ \ldots\subset\CF_{-1}\subset\CF_0\subset\CF_1\subset\ldots$ such that:

(a) $\CF_{k+n}=\CF_k(\bD_0)$ for any $k$;

(b) $ch_1(\CF_k)=k[\bD_0]$ for any $k$: the first Chern classes are
proportional to the fundamental class of $\bD_0$;

(c) $ch_2(\CF_k)=a_i$ for $i\equiv k\pmod{n}$;

(d) $\CF_0$ is locally free at $\bD'$ and trivialized at $\bD':\
\CF_0|_{\bD'}=V\otimes\CO_{\bD'}$;

(e) For $-n\leq k\leq0$ the sheaf $\CF_k$ is locally free at $\bD'$,
and the quotient sheaves $\CF_k/\CF_{-n},\ \CF_0/\CF_k$ (both supported at
$\bD_0=\bC\times \by\subset\CS'$) are both locally free at the point $\bc\times \by$;
moreover, the local sections of $\CF_k|_{\bc\times \bX}$ are those sections
of $\CF_0|_{\bc\times \bX}=V\otimes\CO_\bX$ which take value in $V^k$ at
$\by\in \bX$.

\medskip

We say that a parabolic sheaf $\CF_\bullet$ is {\em locally free} if
$\CF_k$ is locally free for any $k$. Note that this condition implies
that for any $k\leq l\leq k+n$ the quotient sheaf $\CF_l/\CF_k$ is a
locally free sheaf on $\bD_0$, since $\bD_0$ is smooth. Indeed,
$\CF_l/\CF_k$ is a subsheaf in the sheaf
$\CF_{k+n}/\CF_k = \CF_k(\bD_0)/\CF_k = \CF_k(\bD_0)_{|\bD_0}$
which is locally free.

%, and for any $k\leq l\leq k+n$ the
%quotient sheaf $\CF_l/\CF_k$ is a locally free sheaf on $\bD_0$
%(see ~\cite{y} ~3.1).

\subsection{}
\label{CM}
According to ~\cite{hl}, ~\cite{y}, there exists a fine moduli scheme
$\CM^\alpha$ of parabolic sheaves of degree $\alpha$ on $\BA^2$,
and its open subscheme $\oM^\alpha$ which is a fine moduli space of
locally free parabolic sheaves.
We have a natural forgetting morphism $\pi_\alpha:\ \CM^\alpha\to\CA^{a_0},\
(\CF_k)_{k\in\BZ}\mapsto\CF_0$.

\begin{lemma}
\label{yoko}
$\CM^\alpha$ is smooth.
\end{lemma}
{\sl Proof:} 
Let $\Coh$ denote the moduli stack of coherent sheaves on $\bC$ of generic
rank $n$, equipped with a trivialization at $\bc\in\bC$. Let $\FCoh$ denote
the moduli stack of flags $0=\fF_0\subset\fF_1\subset\ldots\subset\fF_n$
of coherent sheaves on $\bC$ with the successive quotients being of generic
rank 1, equipped with a trivialization at $\bc\in\bC$ compatible with the flag.
Both $\Coh$ and $\FCoh$ are smooth, see ~\cite{la}.

We have morphisms $\sr:\ \CA^a\to\Coh,\ 
\CF\mapsto\CF|_{\bD_0}=\CF/\CF(-\bD_0)$, 
and $\fr:\ \CM^\alpha\to\FCoh,\ 
\CF_\bullet\mapsto
\left(\fF_k:=\CF_{k-n}/\CF_0(-\bD_0)\right)_{0\leq k\leq n}$.
Evidently,
$\CM^\alpha$ is the cartesian product of $\CA^{a_0}$ and $\FCoh$ over $\Coh$.

So it remains only to check that $\sr:\ \CA^a\to\Coh$ is smooth.
Both stacks in question being smooth it suffices to show the surjectivity
of the corresponding tangent map. The obstruction to $\sr$ being
a submersion at a point $\CF\in\CA^a$ 
lies in $\Ext^2(\CF,\CF(-\bD'-\bD_0))$. As in the proof of 
Lemma ~\ref{kuz}, it is enough to check $\Hom(\CF,\CF(-\bc\times\bX))=0$.
But exactly this was done in the above cited proof.
\qed

\begin{remark}
\label{dimcm}
We will see in Proposition ~\ref{connect}, Corollary ~\ref{dimensio} 
that $\CM^\alpha$ is connected of dimension
$\dim\CM^\alpha=2|\alpha|:=2\sum_Ia_i$.
\end{remark}

\subsection{}
\label{family}
For a future use we construct a family of regular functions on $\CM^\alpha$
factoring through the projection $\pi_\alpha:\ \CM^\alpha\to\CA^{a_0}$.

Let $\sO$ denote the algebraic variety formed by pairs of lines $(P_1,P_2)$
in the projective plane $\CS''$ such that all the three lines $P_1,P_2,\bD''$
are distinct. Note that $\sO$ is an affine algebraic variety. 
We have a fiber bundle $\sP_1$ (resp. $\sP_2$) over $\sO$ whose fiber over
$(P_1,P_2)$ is $P_1$ (resp. $P_2$). We have a fiber
bundle $\sfp:\ \tO\to\sO$ whose fiber over $(P_1,P_2)$ is 
$P_2^\circ:=P_2-P_2\cap\bD''$
(isomorphic to $\BA^1$). We denote the $a$-th symmetric power of $\tO$ relative
over $\sO$ by $\tO^{(a)}$.

The relative surface $\CS''_\sO:=\CS''\times\sO$ over $\sO$ has two sections
$p_1:=P_1\cap\bD'',\ p_2:=P_2\cap\bD''$, and a relative line $\bD''_\sO$.
Blowing up $p_1,p_2$ and blowing down the proper transform of $\bD''_\sO$
we obtain the relative surface $\CS'_\sO$. Its exceptional divisor
$\bD'_\sO$ is a union of two $\BP^1$-bundles over $\sO$;
in fact, $\CS'_\sO\simeq\sP_1\times_\sO\sP_2$, and $\bD'_\sO=
\sP_1\times_\sO p_2
\bigcup p_1\times_\sO\sP_2$. 

Given a torsion free sheaf $\CF\in\CA^a$ on $\CS''$ trivialized at $\bD''$
we lift it to $\CS''_\sO$, and then apply the relative version of Atiyah's
trick ~\ref{grassmann} to get a torsion free sheaf $\CF'$ on $\CS'_\sO$
trivialized at $\bD'_\sO$. The sheaf $\CF'$ is flat over $\sP_1$, and for
a point $f=(P_1,P_2,c\in P_2)\in\sP_2$ 
its restriction to the fiber $\BP^1_f\simeq P_1$
of $\CS'$ over $f$ is a coherent sheaf on a projective line. 
If $c=P_2\cap\bD''$ then the restriction of $\CF'$ to $\BP^1_f$ is trivialized
by construction: $\CF'|_{\BP^1_f}=V\otimes\CO_{\BP^1_f}$. Since the condition
of triviality is an open condition in the moduli stack of coherent sheaves
on $P_1$, we get a finite subset $D\subset P_2^\circ$ 
such that for $c\not\in D$ the restriction $\CF'|_{\BP^1_f}$ is trivial.

In fact, $D$ is not just a finite subset of $P_2^\circ$ but it carries a
structure of an effective Cartier divisor in $P_2^\circ$. Indeed, the 
restriction of $\CF'$ to the fiber $P_1\times P_2$ of $\CS'_\sO$ over 
$(P_1,P_2)\in\sO$ defines a morphism from $P_2$ to the moduli stack of
coherent sheaves on $P_1$. This stack has a canonical Cartier divisor 
$\Delta_0$ of nontrivial coherent sheaves. We define $D$ as the inverse image
of $\Delta_0$. It is easy to see that $\deg D=a$, and as $(P_1,P_2)$ vary
in $\sO$, these effective divisors form a section $\sD(\CF)$ of $\tO^{(a)}$.   

Let $'\CP^a$ denote the ind-scheme of sections of $\tO^{(a)}$ over $\sO$.
In fact, it is just an infinite-dimensional vector space. 
The above construction defines a morphism $\theta_a:\ \CA^a\to\ '\CP^a$.

\subsection{}
\label{formula}
Using Nakajima's construction of $\CA^a$ as a quiver variety, it is possible
to write down an explicit formula for the map $\theta_a$ above. Recall that
$z$ is our coordinate on $\bC$ which identifies $\bC-\bc$ with $\BA^1$.
Also, $t$ is our coordinate on $\bX$ which identifies $\bX-\bx$ with $\BA^1$.
Thus, $\BA^2=\CS''-\bD''$ is equipped with coordinates $(t,z)$. The variety
$\sO$ is the variety of pairs of nonparallel lines in $\BA^2$. 
If $(P_1^\circ,P_2^\circ)=(\{t=0\},\{z=0\})$, and $\CF\in\CA^\alpha$
is represented by a quadruple $(B_1,B_2,\imath,\jmath)$, then the value
of the section $\theta_a(\CF)$ at the point $(P_1,P_2)$ lies in 
$(\bX-\bx)^{(a)}=\BA^{(a)}$. We will prove that 
${\mathsf{ev}}_{(P_1,P_2)}(\theta_a(\CF))=\Spec(B_2)$
(that is the effective divisor in $\BA^1$ cut out by the equation 
$\det(B_2-t\Id)$), and the morphism $\theta_a$ is equivariant under
the natural action of the group of affine linear transformations of $\BA^2$.

More precisely, given $(P_1,P_2)
\in\sO$ we choose an affine linear transformation $g$ of $\BA^2$ such that
$g(\{t=0\})=P_1^\circ,\ g(\{z=0\})=P_2^\circ$. It also identifies $P_2^\circ$
with $\bX-\bx$ with coordinate $t$.
Let us write $g$ as a linear transformation (with matrix entries 
$g_{11},g_{12},g_{21},g_{22}$) followed by a translation by 
$(g_1,g_2)\in\BA^2$. 
The natural action of affine linear transformations of $\BA^2$ on 
$\CA^a$ in terms of Nakajima's quadruples looks like
$$g(B_1,B_2,\imath,\jmath)=
(g_{11}B_1+g_{12}B_2+g_1\Id,g_{21}B_1+g_{22}B_2+g_2\Id,
(g_{11}g_{22}-g_{12}g_{21})\imath,\jmath)$$
For an operator $B\in\End(W)$ we define
$\Spec_g(B)\in(P_2^\circ)^{(a)}$ as follows. First we consider
an effective divisor on $\BA^2$
cut out by an equation $\det(B-(g_{21}z+g_{22}t-g_2)\Id)$. 
Then we intersect it with $P_2^\circ$.

Now given $\CF\in\CA^a$ we compute 
$\theta_a(\CF,P_1,P_2)\in(P_2^\circ)^{(a)}$ in terms of Nakajima's quiver
data $(B_1,B_2,\imath,\jmath)$ for $\CF$.

The following lemma is borrowed from ~\cite{bgk}.

\begin{lemma}
\label{B12}
$\theta_a(B_1,B_2,\imath,\jmath,P_1,P_2)=\Spec_g(g_{21}B_1+g_{22}B_2)$.
\end{lemma}

{\sl Proof:} Recall that $\CF|_{\BA^2}$ is the middle cohomology of the
monad ~(\ref{monad}). We define the map
$$K_g:\ \Ker(\fb)\to V\otimes\det{}^{-1}
\left(g_{21}B_1+g_{22}B_2-(g_{21}z+g_{22}t-g_2)\Id\right)\CO_{\BA^2}$$ as
follows. It sends $(w_1,w_2,v)\in\Ker(\fb)$ to
$$v-\jmath\left(g_{21}(B_1-z)+g_{22}(B_2-t)+g_2\right)^{-1}
(g_{21}w_1+g_{22}w_2)$$
Since $K_g(\operatorname{Im}(\fa))=0$ we get a well defined map
$$L_g:\ \CF|_{\BA^2}\to V\otimes\det{}^{-1}
\left(g_{21}B_1+g_{22}B_2-(g_{21}z+g_{22}t-g_2)\Id\right)\CO_{\BA^2}$$
It is easy to see that $L_g$ is injective, 
and its image contains $$V\otimes
\det\left(g_{21}B_1+g_{22}B_2-(g_{21}z+g_{22}t-g_2)\Id\right)\CO_{\BA^2}$$
\qed

In particular, in the coordinates $g$, the degree of section $\theta_a(\CF)$ 
as a function of $g$ is less than
or equal to $a$. Let $\CP^a\subset\ '\CP^a$ denote the space of 
sections of $\tO^{(a)}$ of degree less than or equal to $a$. Then our map
$\theta_a:\ \CA^a\to\ '\CP^a$ actually lands into the finite dimensional
subspace $\CP^a$. Note that the morphism $\theta_a$ from $\CA^a$ to the
ind-scheme $'\CP^a$ {\em a priori} lands into a (finite type) subscheme
of $'\CP^a$, and $\CP^a$ is just an explicit estimate of such a subscheme.

Composing $\theta_{a_0}$ with the projection
$\pi_\alpha:\ \CM^\alpha\to\CA^{a_0}$ we get the desired map
$\vartheta_\alpha:\ \CM^\alpha\to\CP^{a_0}$.

When $P_1^\circ=\{t=0\}$, and $P_2^\circ=\{z=0\}$, the fiber of $\tO^{(a)}$
over $(P_1,P_2)$ canonically identifies with $(\bX-\bx)^{(a)}$. Composing
$\theta_a$ with evaluation at $(P_1,P_2)$ we get the map $\eta_a:\
\CA^a\to(\bX-\bx)^{(a)}$. Composing $\eta_{a_0}$ with the projection
$\pi_\alpha:\ \CM^\alpha\to\CA^{a_0}$ we get the map
$\eta_\alpha:\ \CM^\alpha\to(\bX-\bx)^{(a_0)}$.

\subsection{}
\label{mho}
Let $\CF\in\CA^a$ be a torsion free sheaf, let
$\de(\CF)=\underline{s}\in\Sym^d(\BA^2)$
(notations of ~\ref{grassmann}). Then $\CN(\CF)\in
\CA^{a-d}$, and we will compute $\theta_a(\CF)$ in terms
of $\theta_{a-d}(\CN(\CF)),\underline{s}$.

To this end we define a morphism $\mho^d$ from the symmetric power
$\Sym^d(\BA^2)$ to $\CP^d$. Namely, for $\underline{s}\in\Sym^d(\BA^2)$,
and $(P_1,P_2)\in\sO$ we can project $\BA^2$ onto $P_2^\circ$ along 
$P_1^\circ$, and the projection of $\underline{s}$ will be an effective
degree $d$ divisor in $P_2^\circ$. As $(P_1,P_2)$ vary in $\sO$ we get
the desired map $\mho^d:\ \Sym^d(\BA^2)\to\CP^d$.

Also note that we have a natural addition
$\tO^{(k)}\times_\sO\tO^{(l)}\to
\tO^{(k+l)},\ (D_1,D_2)\mapsto D_1+D_2$, which gives rise to the
addition map $\CP^k\times\CP^l\to\CP^{k+l}$. Now we can formulate the
following corollary of Lemma ~\ref{B12} due to Nakajima.
A proof of a similar statement for more general quiver varieties can be
found in ~\cite{n0}, ~3.27; see also ~\cite{l}, ~2.30.

\begin{corollary}
\label{nak}
$\theta_a(\CF)=\mho^d(\underline{s})+\theta_{a-d}(\CN(\CF)).$
\end{corollary}
\qed

\section{Quasimaps into the Kashiwara flag scheme and Uhlenbeck spaces}

\subsection{Based maps and quasimaps}
\label{based maps}
We return to the setup of ~\ref{Pluck}. According to ~\cite{k2},
$H^2(\CB,\BZ)$ is canonically isomorphic to the weight lattice $X:\
\lambda\mapsto c_1(\CL_\lambda)$. The dual lattice $H_2(\CB,\BZ)$ is
canonically isomorphic to the coroot lattice $Y$. We say that a regular
map $\phi:\ \bC\to\CB$ has degree $\alpha\in Y$ if the fundamental class
of $\bC$ in the second homology of $\CB$ equals $\alpha:\ \phi_*[\bC]=\alpha$.
Equivalently, $\deg(\phi)=\alpha$ iff for any $\lambda\in X$ we have
$\deg(\phi^*\CL_\lambda)=\langle\lambda,\alpha\rangle$. Then necessarily
$\alpha\in Y^+=\BN[I]$.

We say that $\phi$ is {\em based} if $\phi(\bc)=B_0$ (notations
of ~\ref{base point}).

According to ~\cite{fg}, for any $\alpha\in Y^+$ there exists a fine
moduli space $\obM^\alpha$ of based maps of degree $\alpha$ from
$(\bC,\bc)$ to $(\CB,B_0)$.
Moreover, it is a smooth connected quasiaffine scheme of
dimension $2|\alpha|$. Let us recall its quasiaffine embedding.
Recall that we have a canonical surjection of vector bundles on $\CB:\
V_\lambda\otimes\CO_\CB\twoheadrightarrow\CL_\lambda$ for any $\lambda\in X^+$.
Dually, we have an embedding of vector bundles:
$\CL_{-\lambda}\hookrightarrow V_\lambda^*\otimes\CO_\CB$. Thus,
$\phi\in\obM^\alpha$ gives rise to a collection of line subbundles
$(\phi^*\CL_{-\lambda}\hookrightarrow
V_\lambda^*\otimes\CO_\bC)_{\lambda\in X^+}$ such that

(a) the fiber of $\phi^*\CL_{-\lambda}$ at $\bc\in \bC$ equals $\ell^0_\lambda
\subset V_\lambda^*$ (notations of ~\ref{base point});

(b) This collection of line subbundles satisfies fiberwise Pl\"ucker equations.

\medskip

Note that $\phi^*\CL_{-\lambda}\simeq\CO_\bC(\langle-\lambda,\alpha\rangle)$.
Hence the datum of $\phi^*\CL_{-\lambda}\hookrightarrow
V_\lambda^*\otimes\CO_\bC\ \Leftrightarrow\
\phi^*\CL_{-\lambda}(\langle\lambda,\alpha\rangle)\hookrightarrow
V_\lambda^*\otimes\CO_\bC(\langle\lambda,\alpha\rangle)$ is equivalent
to the datum of {\em nowhere vanishing} section
$s_\lambda\in\Gamma(C,V_\lambda^*\otimes\CO_\bC(\langle\lambda,\alpha\rangle))$
up to scalar multiplication.

Recall that we have chosen a coordinate $z$ on $\bC-\bc$. Thus $s_\lambda$ is
just a polynomial in $z$ of degree $\langle\lambda,\alpha\rangle$ with
values in $V_\lambda^*$. The condition (a) above means that the scalar
product of $s_\lambda$ with a highest vector of $V_\lambda$
is a (scalar) polynomial of
degree exactly  $\langle\lambda,\alpha\rangle$, and the scalar
product of $s_\lambda$ with any other weight vector of $V_\lambda$
is a (scalar) polynomial of
degree strictly less than $\langle\lambda,\alpha\rangle$. Now we may scale
a constant multiple indeterminacy in the choice of $s_\lambda$ by requirement
that the scalar
product of $s_\lambda$ with {\em the}
highest vector $v_\lambda\in V_\lambda$ is a {\em monic} polynomial of
degree $\langle\lambda,\alpha\rangle$.

All in all, $\obM^\alpha$ is formed by collections of such $V_\lambda^*$-valued
{\em nowhere vanishing} polynomials $s_\lambda$ satisfying Pl\"ucker equations.
If we drop the nowhere vanishing condition, we obtain an affine closure
$\bM^\alpha\supset\obM^\alpha$.

Equivalently, $\bM^\alpha$ is formed by collections of {\em invertible
subsheaves} $(\fL_{-\lambda}\subset
V_\lambda^*\otimes\CO_\bC)_{\lambda\in X^+}$ such that

(a) $\fL_{-\lambda}$ is a line subbundle at $\bc\in \bC$,
and its fiber equals $\ell^0_\lambda
\subset V_\lambda^*$ (notations of ~\ref{base point});

(b) This collection of invertible subsheaves
satisfies fiberwise Pl\"ucker equations;

(c) $\deg(\fL_{-\lambda})=-\langle\lambda,\alpha\rangle$.

\medskip

The points of $\bM^\alpha$ will be called {\em based quasimaps of degree
$\alpha$}.

\subsection{Relative based quasimaps}
\label{rel}
We will need a slight generalization of the above construction.
Let $Q$ be a scheme, let $K$ be a set of indices, and for $k\in K$ let
$W^k$ be a (pro)finite dimensional vector bundle over $Q$.
Let $\BP(W^k)$ be the corresponding projective scheme over $Q$.
Let $R\subset\prod_{k\in K}\BP(W^k)$ (cartesian product over $Q$) be a
closed subscheme. Let $s:\ Q\to R,\ q\mapsto(w^k(q))$ be a section.

A {\em relative based quasimap} $\phi$ from $(\bC,\bc)$ to $(R,s)$ is the
following collection of data:

(a) a point $q\in Q$;

(b) an invertible subsheaf $\fL_k\subset W^k_q\otimes\CO_\bC$ for any $k\in K$

satisfying the following conditions:

(i) for an open subset $U\subset\bC$ the invertible subsheaves $\fL_k\subset
W^k_q\otimes\CO_U$ are {\em line subbundles}, so they give rise to a map
$\phi:\ U\to\prod_{k\in K}\BP(W^k_q)$, and its image is required to lie in $R$;

(ii) we have $\bc\in U$, and $\phi(\bc)=s(q)$.

\medskip

The arguments of the previous subsection show the existence of the
classifying scheme $M(\bC,\bc;R,s)$ for the relative based quasimaps.
Simultaneously we obtain
the open subscheme $\ooM(\bC,\bc;R,s)\subset M(\bC,\bc;R,s)$
classifying the relative based maps (i.e. when $U$ above equals $\bC$).

Note that if $R'\hookrightarrow R$ is a closed subscheme, and the section
$s$ factors through $s:\ Q\to R'\hookrightarrow R$, then
$M(\bC,\bc;R',s)$ is a closed subscheme of $M(\bC,\bc;R,s)$.

\subsection{Quasimaps into Grassmannian}
\label{qg}
Recall the setup of ~\ref{kg}. We have a closed subscheme 
$\CG\subset\BP(V_{\omega_0}^*)$ with the base point $G_0\subset\CG$.
It is well known that $H^2(\CG,\BZ)=\BZ$ is generated by $ch_1(\CL_0)$.
Thus for a positive integer $a$ we have the classifying scheme
$\bA^a\supset\obA^a$ of based quasimaps (resp. maps) of degree $a$ from
$(\bC,\bc)$ to $(\CG,G_0)$.

\subsection{}
\label{identification}
We will construct an identification $\oM^\alpha\equiv\obM^\alpha$.
Let $\CF_\bullet\in\CM^\alpha$
be a parabolic sheaf of degree $\alpha$ on $\BA^2$ (notations of ~\ref{parab}).
Then $\CF_0$ is trivialized at $\bD'$, and in particular, at
$\bC\times \bx\subset\bD':\
\CF_0|_{\bC\times \bx}\stackrel{\sigma}{\longrightarrow}
V\otimes\CO_\bC$.

\begin{lemma}
\label{trivialization}
The trivialization $\sigma$ extends canonically to a
trivialization in the formal neighbourhood of $\bC\times \bx$ in $\CS':\
\CF_0|_{\CS'_{\widehat{\bC\times \bx}}}\stackrel{\varsigma}{\longrightarrow}
V\otimes\CO_{\CS'_{\widehat{\bC\times \bx}}}$.
\end{lemma}

{\sl Proof:} We essentially repeat the arguments of ~\ref{family}.
Let us denote by $\varrho_\bX:\ \CS'=\bC\times \bX\to \bX$
the canonical projection
from $\CS'$ to $\bX$. For any point $y\in \bX$ the fiber
$\varrho_y:=\varrho_\bX^{-1}(y)$ is identified with $\bC$. The sheaf $\CF_0$
is flat over $\bX$, and the restrictions $\CF_0|_{\varrho_y}$ are coherent
sheaves on $\bC$. Thus we obtain a morphism from $\bX$ to the moduli stack
of coherent sheaves on $\bC$. By our assumption this morphism
sends $\bx\in \bX$
to the class of trivial vector bundle on $\bC$. Let $D\subset \bX-\bx$
be the inverse image of the
Cartier divisor of nontrivial coherent sheaves in this moduli stack.
Then our trivialization $\sigma$ extends to a trivialization on
$\bC\times(\bX-D):\
\CF_0|_{\bC\times(\bX-D)}\simeq V\otimes\CO_{\bC\times(\bX-D)}$.
There is a unique choice of such an extension $\varsigma$ such that
$\varsigma|_{\bc\times(\bX-D)}$ coincides with the given trivialization of
$\CF_0|_{\bc\times(\bX-D)}$
(as $\bc\times(\bX-D)\subset\bD'$). Finally, we just
restrict our canonical rational trivialization $\varsigma$ to the formal
neighbourhood of $\bC\times \bx$ in $\CS'$.
\qed

Given a locally free parabolic sheaf $\CF_\bullet\in\oM^\alpha$,
we equip it with the canonical trivialization $\varsigma$ in the formal
neighbourhood of $\bC\times \bx$ in $\CS'$. Consider the canonical projection
$\varrho_\bC$ from $\CS'=\bC\times \bX$ to $\bC$.
For a point $c\in \bC$, restricting
$(\CF_\bullet,\varsigma)$ to the fiber $\varrho_c:=\varrho_\bC^{-1}(c)=\bX$
we obtain a parabolic vector bundle $(\fF_\bullet(c),\tau(c))$ on $\bX$
(notations of ~\ref{P1 moduli}). For $q=\bc$ the corresponding parabolic
vector bundle $(\fF_\bullet(\bc),\tau(\bc))$ gives the point $B_0\in\CB$.
Thus, starting from $\CF_\bullet\in\oM^\alpha$, we have constructed a
based map $\phi$ from $(\bC,\bc)$ to $(\CB,B_0)$. It is easy to see that
$\deg(\phi)=\alpha$, so we have constructed a morphism $\zeta:\ \oM^\alpha\to
\obM^\alpha$.

Conversely, the data of a based map $\phi\in\obM^\alpha$, by the very
definition, consists of a family of
parabolic vector bundles over $\bC$, that is, a locally free parabolic sheaf
$\CF_\bullet$ on $\BA^2$ along with a trivialization of $\CF_0$ in the
formal neighbourhood $\CS'_{\widehat{\bC\times \bx}}$, compatible with a given
trivialization of $\CF_0|_{\bD'}$. To define $\zeta^{-1}(\phi)$ we just
forget the formal trivialization. Thus we have constructed the inverse
isomorphism $\zeta^{-1}:\ \obM^\alpha\to\oM^\alpha$.

\medskip

The same argument establishes an identification $\oCA^a\equiv\obA^a$.

\subsection{Uhlenbeck space}
\label{fA}
Recall that Nakajima defines $\CA^a$ as the moduli space of stable
quadruples $(B_1,B_2,\imath,\jmath)$, see ~\ref{nakajima}. He also defines
the {\em Uhlenbeck space} 
$\fN^a$ as the GIT quotient of the space of quadruples with respect to the
natural GL$(W)$-action. There is a natural proper morphism $\Upsilon_a:\
\CA^a\to\fN^a$ which is a semismall resolution of singularities, 
see ~\cite{n0}. In fact, $\fN^a$ is the affinization of $\CA^a$, that is the
spectrum of the algebra of regular functions on $\CA^a$.

We propose two more definitions of the Uhlenbeck space; we conjecture
that they are both equivalent to Nakajima's definition, see ~\ref{denis}.
 From now on we will identify $\obA^a$ with $\oCA^a$. In particular,
we have an open embedding $\bj:\ \oCA^a\hookrightarrow\bA^a$.
Recall also the morphism $\theta_a:\ \oCA^\alpha\to\CP^a$
defined in ~\ref{formula}. 

We define the Uhlenbeck space $\fA^a$
as the closure of $\oCA^a$ in $\bA^a\times\CP^a$
(with respect to the locally closed embedding $(\bj,\theta_a)$). 

\subsection{}
\label{sA}
The virtue of our second construction of the Uhlenbeck space is that it
carries a natural action of the group of affine linear transformations
of $\BA^2$. Recall the setup of ~\ref{family}. For a projective line 
$P_2$ with a point $p_2=P_2\cap\bD''$ we consider the moduli space
$\CG_{P_2,p_2}$ of $SL_n$-bundles on $P_2$ trivialized in the formal 
neighbourhood of $p_2\in P_2$. It carries an ample line bundle 
$\CL_{P_2,p_2}$ whose fiber at $(\fF,\tau)$ is $\det R\Gamma(P_2,\fF)$
(where we view $\fF$ as a vector bundle of rank $n$ on $P_2$). Also, we
have a base point $G_{P_2,p_2}\in\CG_{P_2,p_2}$, namely, $G_{P_2,p_2}=
(\fF,\tau)$ where $\fF=V\otimes\CO_{P_2}$, and $\tau$ is the tautological
trivialization. Thus we have a relative scheme $\CG_\sO$ over $\sO$ together
with a section $G_\sO$.

Following ~\ref{rel} we define a {\em relative based quasimap of degree $a$} 
from 
$(P_1,p_1)$ to $(\CG_\sO,G_\sO)$ as the following collection of data:

(a) a point $(P_1,P_2)\in\sO$;

(b) an invertible subsheaf $\fL\subset\Gamma^*(\CG_{P_2,p_2},\CL_{P_2,p_2})
\otimes\CO_{P_1}$ of degree $-a$.

satisfying the conditions (i),(ii) of {\em loc. cit.}

\medskip

The arguments of {\em loc. cit.} establish the existence of the classifying
scheme $M^a(P_1,p_1;\CG_\sO,G_\sO)$ for the relative based quasimaps.
We have an evident projection $M^a(P_1,p_1;\CG_\sO,G_\sO)\to\sO$ with a
typical fiber isomorphic to $\bA^a$. All the schemes in question are affine.
Let $\sM^a$ be the ind-scheme of sections of $M^a(P_1,p_1;\CG_\sO,G_\sO)$
over $\sO$. 

Given a locally free sheaf $\CF\in\oCA^a$ we get a locally free sheaf 
$\CF'$ on $\CS'_\sO$ trivialized at $\bD'_\sO$ as in ~\ref{family}.
Applying the relative version of the arguments in ~\ref{identification}
to $\CF'$ we produce from it a section of $M^a(P_1,p_1;\CG_\sO,G_\sO)$
over $\sO$. Thus we construct a morphism $\fj:\ \oCA^a\to\sM^a$.

\medskip

Finally, we define our second version of the Uhlenbeck space $\sA^a$
as the closure of $\oCA^a$ in $\sM^a\times\CP^a$ with respect to the
locally closed embedding $(\fj,\theta_a)$.

Evaluating the sections in $\sM^a$ at the point $(P_1,P_2)=(\{t=0\},\{z=0\})$
we obtain the morphisms $\sM^a\to\bA^a$ and $\Xi_a:\ \sA^a\to\fA^a$.

\subsection{Uhlenbeck flag space}
\label{main char}
 From now on we will identify $\obM^\alpha$ with $\oM^\alpha$. In particular,
we have an open embedding $\bj:\ \oM^\alpha\hookrightarrow\bM^\alpha$.
Recall also the morphism $\vartheta_\alpha:\ \oM^\alpha\to\CP^{a_0}$
defined in ~\ref{formula}.
We are finally able to introduce our main character.

We define\footnote{It was V.Drinfeld 
who noticed that the closure $'\fM^\alpha$
of $\oM^\alpha$ in $\bM^\alpha\times(\bX-\bx)^{(a_0)}$ is a wrong candidate
for the Uhlenbeck flag space; in particular, $'\fM^\alpha$ is not normal in
general.}
the {\em Uhlenbeck flag space} $\fM^\alpha$ as the closure of
$\oM^\alpha$ in $\bM^\alpha\times\CP^{a_0}$
(with respect to the locally
closed embedding $(\bj,\vartheta_\alpha)$).

\begin{proposition}
\label{finite type}
$\fM^\alpha$ is an irreducible affine scheme of finite type, of dimension
$\dim(\fM^\alpha)=2|\alpha|$.
\end{proposition}

{\sl Proof:} According to ~\cite{fg} (alternatively,
see Proposition ~\ref{connect}, Corollary ~\ref{dimensio}), 
$\oM^\alpha$ is smooth, connected
of dimension $2|\alpha|$. Hence it only remains to prove that $\fM^\alpha$
is of finite type, the problem being that $\bM^\alpha$ is not of finite type,
see {\em loc. cit.} Recall the morphism $\eta_\alpha:\ 
\oM^\alpha\to(\bX-\bx)^{(a_0)}$
defined in ~\ref{formula}. 
Let us consider the closure $'\fM^\alpha$
of $\oM^\alpha$ in $\bM^\alpha\times(\bX-\bx)^{(a_0)}$ (with
respect to the locally closed embedding $(\bj,\eta_\alpha)$).
Since $\eta_\alpha$ factors through $\vartheta_\alpha$, we may equivalently
define $\fM^\alpha$ as the closure of $\oM^\alpha$ in $'\fM^\alpha\times
\CP^{a_0}$. Thus it suffices to prove that $'\fM^\alpha$ is of finite type.

Let us denote $a_0$ by $a$ till the end of the proof.
Recall the Beilinson-Drinfeld-Kottwitz ind-scheme $\fB^a$. We have a
closed embedding\footnote{We just mean that the restriction of $(m_a,p_a)$
to a finite type closed subscheme of $\fB^a$ is a closed embedding.}
$(m_a,p_a):\ \fB^a\hookrightarrow\CB\times(\bX-\bx)^{(a)}$
(notations of ~\ref{BDK}). Recall the notion of relative based
(quasi)maps, see ~\ref{rel}. We apply it to the case
$Q=(\bX-\bx)^{(a)}, R=\CB\times(\bX-\bx)^{(a)}$ with the (relative)
Pl\"ucker embedding, and the evident section $s$. Then evidently
$M^\alpha(\bC,\bc;R,s)=\bM^\alpha\times(\bX-\bx)^{(a)}$, and
$\ooM^\alpha(\bC,\bc;R,s)=\obM^\alpha\times(\bX-\bx)^{(a)}$.
%Note in particular that $\ooM^\alpha(\bC,\bc;R,s)$ is dense in
%$M^\alpha(\bC,\bc;R,s)$ since $\obM^\alpha$ is dense in
%$\bM^\alpha$ (see ~\cite{fg}). 
We will use an embedding
$(\Id,\eta_\alpha):\ \obM^\alpha\hookrightarrow\obM^\alpha(\bC,\bc;R,s)$.

\begin{lemma}
\label{Gr+}
There exists a closed subscheme of finite type $\fB^\alpha_+\subset\fB^a$
such that for any based map $\phi\in\obM^\alpha$ the relative based map
$\phi':=(\Id,\eta_\alpha)(\phi)$
factors through $\fB^\alpha_+$.
\end{lemma}

{\sl Proof:} For a fixed $\phi$, the proof of Lemma ~\ref{trivialization}
shows that $\phi'$ factors through
$\fB^a\hookrightarrow\CB\times(\bX-\bx)^{(a)}$, and hence through its
closed subscheme of finite type. We have to choose such a subscheme
uniformly, as $\phi$ varies. Recall that $\obM^\alpha=\oM^\alpha$ is a scheme
of finite type, and we have a natural evaluation morphism ${\mathsf{ev}}:\
\obM^\alpha\times\bC\to\fB^\alpha,\ (\phi,c)\mapsto\phi'(c)$. By definition,
such a morphism into the ind-scheme $\fB^\alpha$ factors through a finite
type scheme $\fB^\alpha_+$. 
\qed

\begin{remark}
Let us give a more concrete description of $\fB^\alpha_+$.
Let $(y_1,\ldots,y_a)\in(\bX-\bx-\by)^{(a)}-\Delta$
be a collection of {\em distinct}
points. Then the fiber of $\fB^a$ over $(y_1,\ldots,y_a)$ equals the
product of $\sB$ (the flag variety of $\fsl_n$) and of
$a$ copies of the affine Grassmannian $\Gr$ of $\fsl_n$,
see ~\cite{g}. Let $\Gr_{\gamma_0}\subset\Gr$ be the closure of 
$SL_n(\CO)$-orbit
numbered by the highest (co)root ${\gamma_0}$ of $\fsl_n$.
Let $\fB^a_1\subset p_a^{-1}((\bX-\bx-\by)^{(a)}-\Delta)$ be a closed
subscheme of finite type whose fiber over
$(y_1,\ldots,y_a)\in(\bX-\bx-\by)^{(a)}-\Delta$ equals
$\sB\times\prod\Gr_{\gamma_0}
\subset\sB\times\prod\Gr$.
Let $\fB^a_2\subset\fB^a$ be the closure of $\fB^a_1$.

The proof of Lemma ~\ref{B12} actually shows that if
$\eta_\alpha(\phi)\in(\bX-\bx-\by)^{(a)}$ then $\phi'$ factors through
$\fB^a_2$.

According to ~\cite{g}, the fiber of $\fB^a$ over $a\cdot\by\in(\bX-\bx)^{(a)}$
equals the affine flag variety $\Fl$ of $\fsl_n$. We define a closed
subscheme of finite type $\Fl_a\subset\Fl$ as the fiber of $\fB^a_2$ over
$a\cdot\by\in(\bX-\bx)^{(a)}$.

We have an embedding
$(\bX-\bx-\by)^{(a-1)}-\Delta\hookrightarrow(\bX-\bx)^{(a)}-\Delta,\
(y_1,\ldots,y_{a-1})\mapsto(\by,y_1,\ldots,y_{a-1})$. According to ~\cite{g},
the fiber of $\fB^a$ over $(\by,y_1,\ldots,y_{a-1})$ equals the product
of $\Fl$ and $a-1$ copies of $\Gr$.
We define $b:=\max(a_i)_{i\in I}$. Let $\fB^a_3\subset
p_a^{-1}((\bX-\bx-\by)^{(a-1)}-\Delta)$ be a closed subscheme of finite
type whose fiber over $(\by,y_1,\ldots,y_{a-1})$ equals
$\Fl_b\times\prod\Gr_{\gamma_0}$. Let $\fB^a_4\subset\fB^a$ be the closure
of $\fB^a_3$.

Finally, we define $\fB^\alpha_+:=\fB^a_2\bigcup\fB^a_4$. Clearly, it
is a closed subscheme of finite type solving our problem.
\end{remark}

\subsection{}
\label{dok-vo}
We return to the proof of Proposition ~\ref{finite type}.
In the setup of ~\ref{rel} we set
$R':=\fB^a_+\stackrel{(m_a,p_a)}{\hookrightarrow}\CB\times(\bX-\bx)^{(a)}=:R$.
Thus we have a closed embedding $M(\bC,\bc;R',s)\hookrightarrow
M(\bC,\bc;R,s)$. Since $R'$ is a scheme of finite type,
$M(\bC,\bc;R',s)$ is a scheme of finite type as well.
By Lemma ~\ref{Gr+}, $\ooM(\bC,\bc;R',s)$ coincides with
the image of $\obM^\alpha$ under the embedding $(\Id,\eta_\alpha)$.
Thus, the closure $'\fM^\alpha$ of $(\Id,\eta_\alpha)(\obM^\alpha)$ in
$M(\bC,\bc;R,s)$ is a closed subscheme in a scheme of finite type
$M(\bC,\bc;R',s)$. Hence $'\fM^\alpha$ is a scheme of finite type itself.
This completes the proof of Proposition ~\ref{finite type}.
\qed

\subsection{}
The same proof as above (using the Beilinson-Drinfeld ind-scheme
$\fG^a$, see ~\ref{kg}, instead of $\fB^\alpha$) shows that both $\sA^a$
and $\fA^a$ are irreducible affine schemes of finite type.

\section{Resolution of singularities $\varpi_\alpha:\ \CM^\alpha\to\fM^\alpha$}

\subsection{Determinant line bundles}
\label{det}
We start with a construction of a morphism $\omega^\alpha:\
\CM^\alpha\to\bM^\alpha$. So let $\CF_\bullet$ be an $S$-point of
$\CM^\alpha$, that is, a parabolic sheaf on $S\times\CS'$ flat over $S$
along with a trivialization $\sigma$
of $\CF_0$ at $S\times\bD'$. We have to construct
the invertible subsheaves
$\fL_{-\omega_i}\subset V_{\omega_i}^*\otimes\CO_{S\times\bC}$ satisfying
the Pl\"ucker equations. Equivalently, we have to construct the generically
surjective maps
$V_{\omega_i}\otimes\CO_{S\times\bC}\to\fL_{\omega_i}=:\fL_i$ defined up to
scalar multiplication. Recall that the fundamental representation
$V_{\omega_i}$ is canonically embedded into the semi-infinite wedge power
$\CQ$ (see ~\ref{wedge}). Hence it suffices to construct the generically
surjective maps $p_i:\ \CQ\otimes\CO_{S\times\bC}\to\fL_i$ defined up to
scalar multiplication.

According to Lemma ~\ref{trivialization}, the trivialization $\sigma$
restricted to $S\times\bC\times\bx$ canonically extends to a trivialization
$\varsigma$ of $\CF_0$ (and hence $\CF_k,\ k\in\BZ$) in the formal
neighbourhood of $S\times\bC\times\bx$ in $S\times\CS'$.

Let us denote by $\CU\subset S\times\CS'$ the open subset
$S\times\CS'-S\times\bC\times\bx$. Let us denote by $\bp:\
S\times\CS'\to S\times\bC$ the natural projection,
and by $\cp:\ \CU\to S\times\bC$ its
restriction to $\CU$. Let us denote by $\oU$ the
intersection of $\CU$ with the formal neighbourhood of $S\times\bC\times\bx$
in $S\times\CS'$ (the ``pointed formal neighbourhood''), and by $\op:\
\oU\to\bC$ the natural projection. Then for any $k\in\BZ$ the trivialization
$\varsigma$ identifies $\op_*(\CF_k|_\oU)$ with
$\bV\otimes\CO_{S\times\bC}$ (notations of ~\ref{wedge}), and
$\sF_k:=\cp_*(\CF_k|_\CU)$ is naturally a discrete lattice in
$\op_*(\CF_k|_\oU)=\bV\otimes\CO_{S\times\bC}$. Recall that
$L_V\otimes\CO_{S\times\bC}:=V\otimes\BC[[t^{-1}]]\otimes\CO_{S\times\bC}$
is a compact lattice in $\bV\otimes\CO_{S\times\bC}$.

\medskip

We define $\fL_i$ as the dual of the determinant line bundle of a natural
{\em Fredholm operator} $(L_V\otimes\CO_{S\times\bC})\oplus\sF_i\to
\bV\otimes\CO_{S\times\bC}$ (notations of ~\cite{bbe}).

\subsection{}
\label{projection}
We still have to construct the generically
surjective maps $p_i:\ \CQ\otimes\CO_{S\times\bC}\to\fL_i$. Recall that
$\CQ$ is a union of finite dimensional subspaces $\CQ^{L_{1,2}}$
(see ~\ref{wedge}). It suffices to construct a compatible system of
maps $p_i^{L_{1,2}}:\ \CQ^{L_{1,2}}\otimes\CO_{S\times\bC}\to\fL_i$.

For small enough compact lattices $L_1,L_2$
(such that $L_2\subset L_1^\perp\subset\bV^*$) we have
$L_1\bigcap\sF_i=0$, and $L_2^\perp+\sF_i=\bV\otimes\CO_{S\times\bC}$
for any $0\leq i\leq n-1$. In effect, by \v{C}ech calculation, this is
equivalent to $R^0\bp_*\CF_i(-N)=0$ and $R^1\bp_*\CF_i(N)=0$ for $N\gg0$.
We define
$\sF_i^{L_2}:=\Ker(L_2^\perp\oplus\sF_i\to\bV\otimes\CO_{S\times\bC})$.
This is a coherent sheaf flat over $S\times\bC$
equipped with a canonical embedding into a vector bundle
$L_2^\perp/L_1\otimes\CO_{S\times\bC}$. Note that $\det(\sF_i^{L_2})=
\fL_i^*\otimes\det(L_2^\perp)$, thus
$\det^*(\sF_i^{L_2})=\fL_i\otimes\det^*(L_2^\perp)$.

The embedding $\sF_i^{L_2}\hookrightarrow L_2^\perp/L_1\otimes\CO_{S\times\bC}$
gives rise to an inverbible subsheaf $\det(\sF_i^{L_2})\subset
\Lambda^\bullet(L_2^\perp/L_1)\otimes\CO_{S\times\bC}$. Dually, we have a
generically surjective morphism
$\Lambda^*(L_2^\perp/L_1)\otimes\CO_{S\times\bC}\to\det^*(\sF_i^{L_2})$,
or equivalently,
$\Lambda^*(L_2^\perp/L_1)\otimes\det(L_2^\perp)\otimes\CO_{S\times\bC}
\to\det^*(\sF_i^{L_2})\otimes\det(L_2^\perp)=\fL_i$.
Now recall that $\Lambda^*(L_2^\perp/L_1)\otimes\det(L_2^\perp)$ is
canonically isomorphic to $\CQ^{L_{1,2}}$ (see ~\ref{wedge}), hence we
have obtained the desired morphism
$p_i^{L_{1,2}}:\ \CQ^{L_{1,2}}\otimes\CO_{S\times\bC}\to\fL_i$.

\subsection{}
\label{define}
Over an open subset $U\subset S\times\bC$ such that
$\CF_\bullet|_{\bp^{-1}(U)}$ is a parabolic vector bundle, the above
construction reduces to $\CQ\otimes\CO_U\twoheadrightarrow\fL_i:=
(\CQ\otimes\CO_U)_{\sF_{W,i}}$ (notations of ~\ref{flag}). Hence
$\fL_{-\omega_i}|_U\subset V_{\omega_i}^*\otimes\CO_U$ satisfy Pl\"ucker
relations, hence $\fL_{-\omega_i}\subset V_{\omega_i}^*\otimes\CO_{S\times\bC}$
satisfy Pl\"ucker relations. Evidently, $U\supset S\times\bc$, and the
fibers of $\fL_{-\omega_i}$ at $S\times\bc$ are as prescribed. So all in all
we have constructed the desired morphism $\omega^\alpha:\ \CM^\alpha\to
\bM^\alpha$. By the same token, we have constructed the morphism
$\upsilon^a:\ \CA^a\to\bA^a$.

\medskip

Recall the morphism $\vartheta_\alpha:\ \CM^\alpha\to\CP^{a_0}$
%(resp. $\eta_\alpha:\ \CM^\alpha\to(\bX-\bx)^{(a_0)}$) 
constructed in ~\ref{formula}.

We define the morphism
$\varpi_\alpha:=(\omega^\alpha,\vartheta_\alpha):\
\CM^\alpha\to\bM^\alpha\times\CP^{a_0}$.
Since $\oM^\alpha$ is dense in $\CM^\alpha$, the morphism $\varpi_\alpha$
factors through $\CM^\alpha\to\fM^\alpha\hookrightarrow
\bM^\alpha\times\CP^{a_0}$. Thus we obtain the
same named morphism $\varpi_\alpha:\ \CM^\alpha\to\fM^\alpha$.

\begin{proposition}
$\varpi_\alpha$ is a proper morphism.
\end{proposition}

{\sl Proof:} Recall that $'\fM^\alpha$ is the closure
of $\oM^\alpha$ in $\bM^\alpha\times(\bX-\bx)^{(a_0)}$ (with
respect to the locally closed embedding $(\bj,\eta_\alpha)$). The morphism
$'\varpi_\alpha:=(\omega^\alpha,\eta_\alpha):\
\CM^\alpha\to\bM^\alpha\times(\bX-\bx)^{(a_0)}$
factors through the morphism $'\varpi_\alpha:\
\CM^\alpha\to\ '\fM^\alpha$, and $'\varpi_\alpha$ factors through
$\CM^\alpha\stackrel{\varpi_\alpha}{\longrightarrow}\fM^\alpha\to\
'\fM^\alpha$, so it suffices to check that $'\varpi_\alpha$
is projective.

Let us consider the moduli scheme $\BM^\alpha\supset\bM^\alpha$ of quasimaps
from $\bC$ to $\CB$ removing the based condition in the definition of
$\bM^\alpha$. Recall that $\varrho_\bX:\ \CS'=\bC\times\bX\to\bX$ is the
canonical projection. For an effective divisor $D\in(\bX-\bx)^{(a_0)}$ we
denote by $\bD_D$ the effective divisor $\varrho_\bX^{-1}(D)$ in $\CS'$.
Let $\TM^\alpha$ be the moduli ind-scheme of the following data:

(a) a divisor $D\in(\bX-\bx)^{(a_0)}$;

(b) a parabolic sheaf $\CF_\bullet$ of degree $\alpha$ on $\CS'$ such
that $V\otimes\CO_{\CS'}(-\infty\cdot\bD_D)\subset\CF_0\subset 
V\otimes\CO_{\CS'}(\infty\cdot\bD_D)$.

\medskip

Note in particular that $\CF_0$ (and hence all the $\CF_k,\ k\in\BZ$)
are trivialized in a Zariski (and hence in the formal) neighbourhood
of $\bC\times\bx\subset\CS'$.
Now the construction of ~\ref{det}--\ref{projection} defines a morphism
$\Omega^\alpha:\ \TM^\alpha\to\BM^\alpha$.
The construction of ~\ref{family} defines a
locally closed
embedding $\CM^\alpha\to\TM^\alpha$, and we have a cartesian diagram
$$
\begin{CD}
\CM^\alpha @>>> \TM^\alpha \\
@VVV            @VVV       \\
\bM^\alpha @>>> \BM^\alpha
\end{CD}
$$
We denote by $\eta':\ \TM^\alpha\to(\bX-\bx)^{(a_0)}$ the tautological
projection.
Note that $\eta_\alpha$ factors through $\CM^\alpha\to\TM^\alpha
\stackrel{\eta'}{\longrightarrow}(\bX-\bx)^{(a_0)}$.

It suffices to prove
that $(\Omega^\alpha,\eta'):\ \TM^\alpha\to\BM^\alpha\times(\bX-\bx)^{(a_0)}$
is ind-projective. Moreover, it is enough to prove that
$\eta':\ \TM^\alpha\to(\bX-\bx)^{(a_0)}$ is ind-projective.
Let us view a ``universal divisor'' $\bD_D$ as a closed subscheme of
$(\bX-\bx)^{(a_0)}\times\CS'$, projective over $(\bX-\bx)^{(a_0)}$.
Then $\TM^\alpha$ is a closed ind-subscheme of a product of certain
inductive limits of
Quot-schemes over $\bD_{k\cdot D+\by},\ k\to\infty$, cf. ~\cite{fk2}, ~p.164.
These Quot-schemes
being projective, $\eta'$ is ind-projective. This completes the proof of the
Proposition.
\qed

\begin{remark}
Lemma ~\ref{B12} shows that instead of the ind-scheme
$\TM^\alpha$ in the above argument, one could use 
the scheme $'\TM^\alpha$ defined as
$\TM^\alpha$ with the condition (b) being replaced by

('b) $V\otimes\CO_{\CS'}(-\bD_D)\subset\CF_0\subset 
V\otimes\CO_{\CS'}(\bD_D)$.
\end{remark}

\subsection{}
\label{denis}
Recall the setup of ~\ref{fA},~\ref{sA}. The relative version over $\sO$
of the construction ~\ref{projection}--\ref{define} defines the proper
morphisms $\epsilon_a:\ \CA^a\to\sA^a,\ \varepsilon_a:\ \CA^a\to\fA^a$
such that $\varepsilon_a=\Xi_a\circ\epsilon_a$. Since $\sA^a$ (resp. $\fA^a$)
is affine, $\epsilon_a$ (resp. $\varepsilon_a$) factors through the
affinization $\Upsilon_a:\ \CA^a\to\fN^a$, and we get the morphisms
$\Psi_a:\ \fN^a\to\sA^a,\ \Phi_a:\ \fN^a\to\fA^a$ such that
$\Phi_a=\Xi_a\circ\Psi_a$.

We conjecture that all the maps $\Psi_a,\Phi_a,\Xi_a$ are isomorphisms.
This can be checked at the level of $\BC$-points. In effect, it is well known
that Nakajima's space $\fN^a$ has a decomposition into locally closed 
pieces $\fN^a=\bigsqcup\limits_{b\leq a}\oCA^{a-b}\times\Sym^b(\BA^2)$.
On the other hand, recall that 
$\bA^a=\bigsqcup\limits_{b\leq a}\oCA^{a-b}\times\Sym^b(\bC-\bc)$. 
Also, recall the embedding $\mho^b:\ \Sym^b(\BA^2)\hookrightarrow\CP^b$,
see ~\ref{mho}. We have $\fA^a\subset\bA^a\times\CP^a$, and the arguments
of ~\ref{tuhl} below show that $\fA^a$ has a decomposition into locally
closed pieces 
$$\fA^a=\bigsqcup_{b\leq a}\left(\oCA^{a-b}\times\Sym^b(\bC-\bc)\right)
\times_{\Sym^b(\bC-\bc)}\Sym^b(\BA^2)=
\bigsqcup_{b\leq a}\oCA^{a-b}\times\Sym^b(\BA^2)$$
Here the embedding 
$\oCA^{a-b}\times\Sym^b(\BA^2)\hookrightarrow\bA^a\times\CP^a$ goes as follows:
$(\CF;\underline{s})\mapsto(\CF,\varrho_\bC(\underline{s});\theta_{a-b}(\CF)+
\mho^b(\underline{s}))$
where we use the natural projection 
$\varrho_\bC:\ \Sym^b(\BA^2)\to\Sym^b(\bC-\bc)$.

\section{$\IC$ stalks}

\subsection{Uhlenbeck stratification}
\label{tuhl}
We refine the decomposition of the Uhlenbeck space $\fA^a$ into locally
closed pieces described in ~\ref{denis}. 
We define the {\em diagonal stratification} of $\Sym^b(\BA^2)$ as follows. 
For a positive integer $b$ we denote by $\fP(b)$ the set of partitions
of $b$ (in the traditional meaning). For 
$\fP=(b_1\geq b_2\geq\ldots\geq b_m>0)\in
\fP(b)$ the corresponding stratum $\Sym^b(\BA^2)_\fP$ of $\Sym^b(\BA^2)$
is formed by configurations which can be subdivided into $m$ groups of
points, the $r$-th group containing $b_r$ points; all the points in one
group equal to each other, the different groups being disjoint. We have
$\Sym^b(\BA^2)=\bigsqcup\limits_{\fP\in\fP(b)}\Sym^b(\BA^2)_\fP$.

In the setup of ~\ref{denis}, for $b\leq a,\ \fP\in\fP(b)$, 
we define a locally closed subscheme
$\fA^a_{a-b,\fP}=\oCA^{a-b}\times\Sym^b(\BA^2)_\fP\subset
\oCA^{a-b}\times\Sym^b(\BA^2)\subset\bA^a\times\CP^a$.

In order to show 
$\fA^a=\bigsqcup\limits^{\fP\in\fP(b)}_{b\leq a}\fA^a_{a-b,\fP}$ we describe
the inverse image $\varepsilon^{-1}(\fA^a_{a-b,\fP})\subset\CA^a$.
First we define the {\em saturation} and {\em defect} of a based quasimap
$\phi=(\fL_0\subset V^*_{\omega_0}\otimes\CO_\bC)\in\bA^a$. Namely,
the saturation $\CN(\phi)\in\obA^d$ is the based map
$(\widetilde\fL_0\subset 
V^*_{\omega_0}\otimes\CO_\bC)$ where
the line subbundle $\widetilde\fL_0$ is the saturation
of the invertible subsheaf $\fL_0$. The quotient
$\widetilde\fL_0/\fL_0$ is a torsion sheaf on $\bC-\bc$
of length $b$ supported at a finite subset $S$, and we define the
defect $\de(\phi)\in(\bC-\bc)^{(b)}$ as
$\sum_{s\in S}\length_s(\widetilde\fL_0/\fL_0)\cdot s$.
Note that necessarily $a=d+b$.

Recall that the saturation and defect of a torsion free sheaf $\CF\in\CA^a$
were defined in ~\ref{grassmann}. Recall the morphism $\varepsilon_a=
(\upsilon^a,\theta_a):\ \CA^a\to\fA^a\subset\bA^a\times\CP^a$, and the
natural projection $\varrho_\bC:\ \Sym^b(\BA^2)\to(\bC-\bc)^{(b)}$.
We have the following

\begin{lemma}
Suppose the defect $\de(\CF)$ of a torsion free sheaf $\CF\in\CA^a$ has
degree $b\leq a$, so that the saturation $\CN(\CF)$ lies in $\oCA^{a-b}$.
Then the quasimap $\phi:=\upsilon^a(\CF)\in\bA^a$ has saturation
$\CN(\phi)=\upsilon^{a-b}(\CN(\CF))$, 
and defect $\de(\phi)=\varrho_\bC(\de(\CF))\in(\bC-\bc)^{(b)}$. 
%Moreover, $\theta_a(\CF)=\mho^b(\de(\CF))+\theta_{a-b}(\CN(\CF))$.
\end{lemma}

{\sl Proof:} Let $\phi=(\fL_0\subset V^*_{\omega_0}\otimes\CO_\bC)$, and
$\phi':=\upsilon^{a-b}(\CN(\CF))=
(\fL'_0\subset V^*_{\omega_0}\otimes\CO_\bC)$. Since $\CF\subset\CN(\CF)$,
by the construction of $\upsilon$, we see that $\fL_0\subset\fL'_0$.
Moreover, since $\CN(\CF)\in\oCA^{a-b}=\obA^{a-b}$, we see that $\fL'_0$
is saturated; hence $\fL'_0$ is the saturation of $\fL_0$. Since $\CF$
and $\CN(\CF)$ coincide off the support of $\de(\CF)$, by the construction
of $\upsilon$, $\fL_0$ and $\fL'_0$ coincide off the support of
$\varrho_\bC(\de(\CF))$. Finally, since $\length(\fL'_0/\fL_0)=b$, 
and the definition of $\upsilon$ is local over $\bC$, we conclude that
$\de(\phi)=\varrho_\bC(\de(\CF))$.
\qed

\medskip

It follows immediately from the above Lemma and Corollary ~\ref{nak} that 
$\CA^a_{a-b,\fP}:=\varepsilon^{-1}(\fA^a_{a-b,\fP})\subset\CA^a$ is formed
by all the torsion free sheaves $\CF\in\CA^a$ such that $\CN(\CF)$ lies in
$\oCA^{a-b}$, and $\de(\CF)$ lies in $\Sym^b(\BA^2)_\fP$. Clearly, $\CA^a$
is covered by the above locally closed pieces: 
$\CA^a=\bigsqcup\limits^{\fP\in\fP(b)}_{b\leq a}\CA^a_{a-b,\fP}$. Hence
$\fA^a$ being the image of $\CA^a$ under the proper morphism $\varepsilon_a$,
is covered by the Uhlenbeck strata:
$\fA^a=\bigsqcup\limits^{\fP\in\fP(b)}_{b\leq a}\fA^a_{a-b,\fP}$.

\subsection{Saturation and defect}
\label{norm}
Our next goal is to describe a similar stratification of the Uhlenbeck flag
space.

Given a parabolic sheaf $\CF_\bullet\in\CM^\alpha$ we define its
{\em saturation} $\CN(\CF)_\bullet$
as the parabolic sheaf formed by saturations of all components of
the parabolic sheaf $\CF_\bullet$:
$$
\CN(\CF)_k:=\CN(\CF_k).
$$
Clearly, $\CN(\CF)_\bullet$ is indeed a parabolic sheaf. Moreover,
it is evidently locally free.

%First of all, for $k\in\BZ$ we set
%$\CN'(\CF)_k:=\CN(\CF_k)$ (notations of ~\ref{grassmann}). Evidently,
%$\CN'(\CF)_{k+n}=\CN'(\CF)_k(\bD_0)$, so we get a periodic flag of
%vector bundles. Now for $-n\leq k\leq0$ 
%the restriction $\CN'(\CF)_k|_{\bD_0}$
%maps naturally to $\CN'(\CF)_0|_{\bD_0}$. Let us denote the image of this
%morphism by $'\CF_k$. This is a locally free subsheaf in
%$\CN'(\CF)_0|_{\bD_0}$.
%So we have a
%saturation $\widetilde{'\CF_k}\subset\CN'(\CF)_0|_{\bD_0}$ (the minimal
%vector subbundle in $\CN'(\CF)_0|_{\bD_0}$ containing $'\CF_k$).
%Finally, $\CN(\CF)_k\subset\CN'(\CF)_0$ is defined as a subsheaf whose
%local sections are those sections of $\CN'(\CF)_0$ which take value in
%$\widetilde{'\CF_k}\subset\CN'(\CF)_0|_{\bD_0}$ at $\bD_0$. In particular,
%$\CN(\CF)_0=\CN'(\CF)_0=\CN(\CF_0)$, and $\CN(\CF)_{-n}=\CN(\CF)_0(-\bD_0)$.
%Thus, the flag $\CN(\CF)_{-n}\subset\CN(\CF)_{-n+1}\subset\ldots\subset
%\CN(\CF)_0$ can be extended to the desired locally free parabolic sheaf
%$\CN(\CF)_\bullet$.

By construction, $\CF_k\subset\CN(\CF)_k$ for any $k\in\BZ$.

\medskip

Given $\beta=\sum_ib_ii\in\BN[I]$ we define $\BA^\beta=(\bC-\bc)^\beta$ as
the product $\prod_i(\bC-\bc)^{(b_i)}$. This is the moduli space of effective
$I$-colored divisors of degree $\beta$ on $\BA^1=\bC-\bc$.
We define $\oA^2\subset\BA^2$ as $\BA^2-\bD_0$. Note that $\bD_0\cap\BA^2=
\bC-\bc$.

For a parabolic sheaf $\CF_\bullet\subset\CM^\alpha$ we define its {\em defect}
$\de(\CF_\bullet)=
\de^\dagger(\CF_\bullet)+\de^\circ(\CF_\bullet)=
\sum_{i\in I}\de_i(\CF_\bullet)+\de^\circ(\CF_\bullet)
\in\BA^\beta\times\Sym^d(\oA^2)$ (for some $\beta,d$)
as follows. Note that for any $k\in\BZ$ the quotient sheaf
$\CN(\CF)_k/\CF_k|_{\oA^2}$ is a torsion sheaf of finite length $d$ on
$\oA^2$ independent of $k$. In particular, it is supported on a finite
subset $S\subset\oA^2$. So $\de^\circ(\CF_\bullet):=\sum_{s\in S}
\length_s(\CN(\CF)_k/\CF_k)\cdot s\in\Sym^d(\oA^2)$.

Now for $i\in I$ the quotient $\CN(\CF)_i/\CF_i$ is a torsion sheaf
of finite length $d+b_i$. In particular, it is supported on a finite subset
$S_i\subset\BA^2$. So $\de_i(\CF_\bullet):=
\sum_{s\in S_i\cap\bD_0}\length_s(\CN(\CF)_i/\CF_i)\cdot s$, and
$\de^\dagger(\CF_\bullet):=
\sum_i\sum_{s\in S_i\cap\bD_0}\length_s(\CN(\CF)_i/\CF_i)\cdot s\in
\prod_i(\bC-\bc)^{(b_i)}$.

Finally, we define $\de^\ddagger(\CF_\bullet):=
\sum_{s\in S_0}\length_s(\CN(\CF)_0/\CF_0)\cdot s\in\Sym^{d+b_0}(\BA^2)$.

\medskip

Recall that the {\em imaginary coroot} of $\hsl_n$ is $\delta_0:=
\sum_i\alpha_i$. Note that if $\CF_\bullet\subset\CM^\alpha$, and
$\CN(\CF)_\bullet\in\oM^\gamma$, and $\de(\CF_\bullet)\in
\BA^\beta\times\Sym^d(\oA^2)$, then $\alpha=\gamma+\beta+d\delta_0$.

\subsection{Partitions}
\label{part}
We refer the reader to ~\cite{fk2} ~2.2
for a general terminology on partitions.
So for $\alpha\in\BN[I]$ we have the set of {\em usual partitions}
$\Gamma(\alpha)$. The definition of Kostant partitions $\fK(\alpha)$ of
{\em loc. cit.} makes reference to the set $R^+$ of {\em $\psi$-roots},
or in other words, ``positive roots of $\hgl_n$''. To avoid a confusion
we will denote the set of $\psi$-roots by $R^+_{\hgl_n}$, and the
corresponding set of {\em $\hgl_n$-Kostant partitions} by
$\fK_{\hgl_n}(\alpha)$. We define a subset $R^+_{\hsl_n}\subset R^+_{\hgl_n}$
whose complement consists of $\psi$-roots $(0,kn),\ k=1,2,\ldots$
(notations of {\em loc. cit.} ~2.1). The set $R^+_{\hsl_n}$ together
with the dimension function restricted from $R^+_{\hgl_n}$ gives rise to
the set of {\em $\hsl_n$-Kostant partitions}
$\fK_{\hsl_n}(\alpha)$.

The number of summands in a partition $?$ is denoted by $\sK(?)$.

\subsection{Configurations}
\label{config}
We define the {\em diagonal stratification} of $\BA^\alpha=(\bC-\bc)^\alpha$.

Recall that the {\em multisubsets} of a set $S$ are defined as elements of some
symmetric power $S^{(m)}$, and we denote the image of
$(s_1,\dots,s_m)\in S^m$ by $\{\{s_1,\dots,s_m\}\}$.
In particular, the set of usual partitions $\Gamma(\alpha)$ is formed by
all the multisubsets $\Gamma=\{\{\gamma_1,\dots,\gamma_m\}\}$ of $\BN[I]-\{0\}$
such that $\sum_{r=1}^m\gamma_r=\alpha$.

For $\Gamma\in\Gamma(\alpha)$ the corresponding stratum $\BA^\alpha_\Gamma$
is defined as follows. It is formed by configurations which can be
subdivided into $m$ groups of points, the $r$-th group containing $\gamma_r$
points; all the points in one group equal to each other, the different
groups being disjoint. We have $\dim(\BA^\alpha_\Gamma)=\sK(\Gamma)$.
For example, the main diagonal in $\BA^\alpha$
is the closed stratum given by partition $\alpha=\alpha$, while the complement
to all diagonals in $\BA^\alpha$ is the open stratum given by partition
$\alpha=\sum_{i\in I}(\underbrace{i+i+\ldots+i}_{a_i\ \operatorname{ times}})$.
Evidently,
$\BA^\alpha=\bigsqcup\limits_{\Gamma\in\Gamma(\alpha)}\BA^\alpha_\Gamma$.

\medskip

Similarly, we define the {\em diagonal stratification} of $\Sym^d(\oA^2)$.
We have $\Sym^d(\oA^2)=\bigsqcup\limits_{\fP\in\fP(d)}\Sym^d(\oA^2)_\fP$.
Also, in evident notations,
$\Sym^d(\BA^2)=\bigsqcup\limits_{\fP\in\fP(d)}\Sym^d(\BA^2)_\fP$.
For a partition $\fP\in\fP(d)$ let $\Sym^d(\BA^2)_\ofP$ denote the closure
of the stratum $\Sym^d(\BA^2)_\fP$. If $\fP=\{\{k_1\cdot d_1,\ldots,k_m\cdot
d_m\}\}$ for some $0<d_1<\ldots<d_m$, then we have an evident finite
morphism $N^\fP:\ \prod\limits_{l=1}^m\Sym^{k_l}(\BA^2)\to\Sym^d(\BA^2)_\ofP$.
The morphism $N^\fP$ is generically one-to-one; moreover,
$\prod\limits_{l=1}^m\Sym^{k_l}(\BA^2)$ is just the normalization
$\Sym^d(\BA^2)_\tfP$ of the stratum closure $\Sym^d(\BA^2)_\ofP$.
We have $\dim(\Sym^d(\BA^2)_\ofP)=\dim(\Sym^d(\BA^2)_\tfP)=2\sK(\fP)$.

\subsection{Defect stratification}
\label{strat}
We define the {\em defect stratification} of $\CM^\alpha$. Recall the setup
of ~\ref{norm}.
For a decomposition $\alpha=\gamma+\beta+d\delta_0$, and a partition
$\Gamma\in\Gamma(\beta)$, and a partition $\fP\in\fP(d)$, the stratum
$\CM^\alpha_{\gamma,\Gamma,\fP}$ is formed by all $\CF_\bullet\in\CM^\alpha$
such that $\de^\dagger(\CF_\bullet)\in\BA^\beta_\Gamma$, and
$\de^\circ(\CF_\bullet)\in\Sym^d(\oA^2)_\fP$. We have
$\CM^\alpha=\bigsqcup\limits_{\alpha=
\gamma+\beta+d\delta_0}^{\Gamma\in\Gamma(\beta),\fP\in\fP(d)}
\CM^\alpha_{\gamma,\Gamma,\fP}$.

We have an evident projection $\CN:\ \CM^\alpha_{\gamma,\Gamma,\fP}\to
\oM^\gamma$. Also we have a morphism $\de^\ddagger:\
\CM^\alpha_{\gamma,\Gamma,\fP}\to\Sym^{d+b_0}(\BA^2)$.

\medskip

We define the {\em defect stratification} of $\bM^\alpha$ as follows.
For a based
quasimap $\bM^\alpha\ni\phi=
(\fL_{-\omega_i}\subset V^*_{\omega_i}\otimes\CO_\bC)_{i\in I}$
we define its {\em saturation} $\CN(\phi)\in\obM^\gamma$ as the collection
$(\widetilde\fL_{-\omega_i}\subset 
V^*_{\omega_i}\otimes\CO_\bC)_{i\in I}$ where
the line subbundle $\widetilde\fL_{-\omega_i}$ is the saturation
of the invertible subsheaf $\fL_{-\omega_i}$. The quotient
$\widetilde\fL_{-\omega_i}/\fL_{-\omega_i}$ is a torsion sheaf on $\bC-\bc$
of length $b_i$ supported at a finite subset $S_i$, and we define the
{\em defect} $\de(\phi)=\sum_{i\in I}\de_i(\phi)\in\BA^\beta$ as
$\sum_{i\in I}\sum_{s\in S_i}
\length_s(\widetilde\fL_{-\omega_i}/\fL_{-\omega_i})\cdot s$.
Note that necessarily $\alpha=\gamma+\beta$.

Finally, for a decomposition $\alpha=\gamma+\beta$, and a partition
$\Gamma\in\Gamma(\beta)$, the stratum $\bM^\alpha_{\gamma,\Gamma}$ is
formed by all $\phi\in\bM^\alpha$ such that $\de(\phi)\in\BA^\beta_\Gamma$.
Note that $\bM^\alpha_{\gamma,\Gamma}\simeq\obM^\gamma\times\BA^\beta_\Gamma$,
and $\bM^a=\bigsqcup\limits_{\alpha=\gamma+\beta}^{\Gamma\in\Gamma(\beta)}
\bM^\alpha_{\gamma,\Gamma}$ (cf. ~\cite{fg}).

\medskip

Recall that $\varrho_\bC:\ \CS'=\bC\times\bX\to\bC$ is the canonical
projection from $\CS'$ to $\bC$, and $\varrho_\bX:\ \CS'=\bC\times\bX\to\bX$
is the canonical projection from $\CS'$ to $\bX$. Recall also the setup
of ~\ref{formula}, ~\ref{mho}.

\begin{lemma}
\label{comput}
For a parabolic sheaf
$\CF_\bullet\in\CM^\alpha_{\gamma,\Gamma,\fP}$ the quasimap
$\phi:=\omega^\alpha(\CF_\bullet)\in\bM^\alpha$ has saturation
$\CN(\phi)=\omega^\gamma(\CN(\CF)_\bullet)$, and defect
$\de(\phi)=\de^\dagger(\CF_\bullet)+
\delta_0\cdot\varrho_\bC(\de^\circ(\CF_\bullet))$.
%Moreover, $\eta_\alpha(\CF_\bullet)=\varrho_\bX(\de^\ddagger(\CF_\bullet))+
%\eta_\gamma(\CN(\CF)_\bullet)$. 
Furthermore, $\vartheta_\alpha(\CF_\bullet)=
\mho^{d+b_0}(\de^\ddagger(\CF_\bullet))+\vartheta_\gamma(\CN(\CF)_\bullet)$.
\end{lemma}

{\sl Proof:} Clear from definitions and Corollary ~\ref{nak}.
\qed

\subsection{Uhlenbeck flag stratification}
\label{uhl}
For a decomposition $\alpha=\gamma+\beta+d\delta_0$, and a partition
$\Gamma\in\Gamma(\beta)$, and a partition $\fP\in\fP(d)$, we define
a constructible subset
$\fM^\alpha_{\gamma,\Gamma,\fP}\subset\fM^\alpha$ as
$\varpi_\alpha(\CM^\alpha_{\gamma,\Gamma,\fP})$.
Lemma ~\ref{comput} implies that $\CM^\alpha_{\gamma,\Gamma,\fP}=
\varpi_\alpha^{-1}(\fM^\alpha_{\gamma,\Gamma,\fP})$. It follows that
$\fM^\alpha_{\gamma,\Gamma,\fP}$ is a locally closed subscheme of
$\fM^\alpha$, and $\varpi_\alpha:\ \CM^\alpha_{\gamma,\Gamma,\fP}\to
\fM^\alpha_{\gamma,\Gamma,\fP}$ is proper. Moreover, one can see
easily that $\fM^\alpha_{\gamma,\Gamma,\fP}$ is smooth.

Thus we have the Uhlenbeck flag stratification
$\fM^\alpha=\bigsqcup\limits_{\alpha=
\gamma+\beta+d\delta_0}^{\Gamma\in\Gamma(\beta),\fP\in\fP(d)}
\fM^\alpha_{\gamma,\Gamma,\fP}$.

Lemma ~\ref{comput} implies that $\fM^\alpha_{\gamma,\Gamma,\fP}\simeq
\oM^\gamma\times\BA^\beta_\Gamma\times\Sym^d(\oA^2)_\fP$. Let
$\fM^\alpha_{\overline{\gamma,\Gamma,\fP}}$ denote the closure of the
stratum $\fM^\alpha_{\gamma,\Gamma,\fP}$.
Furthermore, if $\beta=0$ (hence $\Gamma=\emptyset$), the above isomorphism
extends to the finite morphism
$\fM^\gamma\times\Sym^d(\BA^2)_\ofP\to
\fM^\alpha_{\overline{\gamma,\emptyset,\fP}}$.
Composing it with the normalization
morphism $N^\fP$ from ~\ref{config} we obtain the finite morphism
$N^{\gamma,\fP}:\ \fM^\gamma\times\Sym^d(\BA^2)_\tfP\to
\fM^\alpha_{\overline{\gamma,\emptyset,\fP}}$ which is generically one-to-one.

\begin{theorem}
\label{semismall}
$\varpi_\alpha:\ \CM^\alpha\to\fM^\alpha$ is semismall. The relevant strata
in $\fM^\alpha$ are the ones with $\beta=0,\ \Gamma=\emptyset$, i.e.
$\fM^\alpha_{\gamma,\emptyset,\fP}$ for $\alpha=\gamma+d\delta_0,\
\fP\in\fP(d)$. The fibers of $\varpi_\alpha$ over the relevant strata are
irreducible.
\end{theorem}

{\sl Proof:} Recall that $\dim(\fM^\alpha)=\dim(\CM^\alpha)=2|\alpha|$.
It follows from the above discussion that
$\dim(\fM^\alpha_{\gamma,\Gamma,\fP})=2|\gamma|+\sK(\Gamma)+2\sK(\fP)$.
The fibers of $\varpi_\alpha:\ \CM^\alpha_{\gamma,\Gamma,\fP}\to
\fM^\alpha_{\gamma,\Gamma,\fP}$
were computed in ~\cite{b}, ~\cite{fk2}, ~\cite{nk}. Namely, consider
a locally free parabolic sheaf $\CF_\bullet\in\oM^\gamma$,
and $\beta'\in\BN[I]$. Recall the
projective variety $K_{\beta'}(\CF_\bullet)$ introduced in ~\cite{fk2}, ~3.1.3.
For a point $c\in\bC-\bc$ let
$\Fib_{\beta',c}(\CF_\bullet)\subset K_{\beta'}(\CF_\bullet)$ be a closed
subvariety given by the condition that the quotient sheaf $T_\bullet$
(see {\em loc. cit.}) is concentrated at $c\in\bC$. For a point
$s\in\oA^2$ and $d'\in\BN$
let $\Fib_{d',s}(\CF_0)$ be the projective variety classifying
all the torsion free subsheaves $\CF'\subset\CF_0$ with $\de(\CF')=d'\cdot s$
(see the ~Appendix of ~\cite{b}).

Let $\fP=\{\{k_1\cdot d_1,\ldots,k_m\cdot d_m\}\}$, and
$\Gamma=\{\{n_1\cdot\beta_1,\ldots,n_g\cdot\beta_g\}\}$ for distinct
$\beta_1,\ldots,\beta_g$. Let $\underline{s}=(d_1\cdot s_1^1,\ldots
d_1\cdot s_1^{k_1},\ldots,d_m\cdot s_m^1,\ldots,d_m\cdot s_m^{k_m})\in
\Sym^d(\oA^2)_\Gamma$, and $\underline{c}=((\beta_1\cdot c_1^1,\ldots
\beta_1\cdot c_1^{n_1},\ldots,\beta_g\cdot c_g^1,\ldots,
\beta_g\cdot c_g^{n_g})\in\BA^\beta_\Gamma$, and $\phi\in\obM^\gamma$
correspond to $\CF_\bullet$ in $\oM^\gamma$.
Then the reduced fiber $\varpi_\alpha^{-1}(\phi,\underline{c},\underline{s})$
is isomorphic to
$$\prod_{l=1}^m(\Fib_{\beta_l,c_l^1}(\CF_\bullet)
\times\ldots\times\Fib_{\beta_l,c_l^{n_l}}(\CF_\bullet))
\times\prod_{l=1}^g(\Fib_{d_l,s_l^1}(\CF_0)
\times\ldots\times\Fib_{d_l,s_l^{k_l}}(\CF_0))$$
Now according to the ~Appendix of ~\cite{b}, the variety
$\Fib_{d',s}(\CF_0)$ is irreducible of dimension $nd'-1$. And according
to ~\cite{fk2}, ~Theorem ~1, $\Fib_{\beta',c}(\CF_\bullet)$ is a union of
various irreducible components of dimension smaller than or equal to
$|\beta'|-1$.

Now a routine arithmetical check completes the proof of the Theorem.
\qed

\begin{corollary}
\label{decomp}
$\varpi_{\alpha*}\uC[2|\alpha|]=
\bigoplus\limits_{\alpha=\gamma+d\delta_0}^{\fP\in\fP(d)}
\IC(\fM^\alpha_{\overline{\gamma,\emptyset,\fP}})$.
\end{corollary}
\qed

\subsection{Symmetric algebras}
\label{direct}
We compute the stalks of $\varpi_{\alpha*}\uC[2|\alpha|]$.
To this end we have to know
the cohomology of fibers of $\varpi_\alpha$. A cellular decomposition
of $\Fib_{\beta',c}(\CF_\bullet)$ is constructed, and the dimensions of the
cells are computed in ~\cite{fk2}, ~p. ~165. The cohomology of
$\Fib_{d',s}(\CF_0)$ (equal to $H^\bullet(\CA(n,d'))$) are computed
in ~\cite{nk}, ~Theorem ~2.3. To arrange the cited information into a neat
form we need some linear algebraic preliminaries.

Recall the {\em Hall algebra} $\bH$ of the category of nilpotent
representations of the cyclic quiver $\tA_{n-1}$, see e.g. ~\cite{fk2}, ~1.3.
It is naturally $\BN[I]$-graded by the dimension of representation
$\bH=\oplus_{\beta\in\BN[I]}\bH_\beta$, see {\em loc. cit.} It has also a
natural filtration $F^0\bH\subset F^1\bH\subset\ldots$, namely, we say that
a class $[M]$ of a nilpotent representation $M$ lies in $F^k\bH$ if
$M$ is a direct sum of $k$ indecomposable representations. In particular,
$F^1\bH/F^0\bH=\BC[R^+_{\hgl_n}]$ (notations of ~\ref{part}).
We denote by $\bH^\bullet$ the associated
graded algebra $\gr_{F^\bullet}\bH$. It is canonically isomorphic to
$\Sym^\bullet(\BC[R^+_{\hgl_n}])$. We have
$\bH^\bullet=\oplus_{\beta\in\BN[I]}\bH^\bullet_\beta$.

Let $\hn_+$ be a subalgebra of $\hsl_n$ generated by $e_i,\ i\in I$.
Choosing a root basis\footnote{To avoid a confusion between roots and
coroots, we should have worked with the {\em Langlands dual} Lie algebra
$\hsl{}_n^\vee$. We prefer to use a coincidence $\hsl_n\simeq\hsl{}_n^\vee$
to save notations at this moment. To clear up things, the interested reader
is referred to ~\ref{speculation}.}
we identify it with $\BC[R^+_{\hsl_n}]$.
Thus we have $\Sym^\bullet(\hn_+)=\Sym^\bullet(\BC[R^+_{\hsl_n}])\subset
\Sym^\bullet(\BC[R^+_{\hgl_n}])=\bH^\bullet$. Also, we have a natural
grading $\Sym^\bullet(\hn_+)=\oplus_{\beta\in\BN[I]}\Sym^\bullet(\hn_+)_\beta$.

\medskip

We also define a bigraded space
$\fu(\hgl_n)=\bigoplus\limits_{d=1,2,\ldots}^{1\leq k\leq n}\fu(\hgl_n)^k_d$
where $\fu(\hgl_n)^k_d$ is 1-dimensional $\BC$-vector space.
We define
$\fu(\hsl_n)=\bigoplus\limits_{d=1,2,\ldots}^{2\leq k\leq n}\fu(\hgl_n)^k_d
\subset\fu(\hgl_n)$. Thus, the symmetric algebras
$\Sym(\fu(\hgl_n)),\ \Sym(\fu(\hsl_n))$ are also bigraded:
$\Sym(\fu(\hgl_n))=\bigoplus\limits_{d\in\BN}^{k\in\BN}\Sym(\fu(\hgl_n))^k_d,\
\Sym(\fu(\hsl_n))=\bigoplus\limits_{d\in\BN}^{k\in\BN}\Sym(\fu(\hsl_n))^k_d$.

\medskip

Recall the notations of the proof of Theorem ~\ref{semismall}.

\begin{proposition}
\label{stalk}
The stalk of $\varpi_{\alpha*}\uC[2|\alpha|]$ at a point
$(\phi,\underline{c},\underline{s})\in\fM^\alpha_{\gamma,\Gamma,\fP}$
is isomorphic to
$$\bigotimes_{l=1}^m\left(\oplus_{r\in\BN}\bH^r_{\beta_l}[2r]\right)^{\otimes
n_l}\otimes\bigotimes_{l=1}^g\left(\oplus_{r\in\BN}\Sym(\fu(\hgl_n))^r_{d_l}
[2r]\right)^{\otimes k_l}[2|\gamma|]$$
\end{proposition}

{\sl Proof:} Follows immediately from the proof of Theorem ~\ref{semismall},
and the above cited results of ~\cite{nk}, ~\cite{fk2}.
\qed

\begin{theorem}
\label{IC}
The stalk of $\IC(\fM^\alpha)$ at a point
$(\phi,\underline{c},\underline{s})\in\fM^\alpha_{\gamma,\Gamma,\fP}$
is isomorphic to
$$\bigotimes_{l=1}^m\left(\oplus_{r\in\BN}
\Sym^r(\hn_+)_{\beta_l}[2r]\right)^{\otimes
n_l}\otimes\bigotimes_{l=1}^g\left(\oplus_{r\in\BN}\Sym(\fu(\hsl_n))^r_{d_l}
[2r]\right)^{\otimes k_l}[2|\gamma|]$$
\end{theorem}

{\sl Proof:} Recall the finite, generically one-to-one morphism
$N^{\gamma,\fP}:\ \fM^\gamma\times\Sym^d(\BA^2)_\tfP\to
\fM^\alpha_{\overline{\gamma,\emptyset,\fP}}$ introduced in ~\ref{uhl}.
It is well known that $\Sym^d(\BA^2)_\tfP$ is rationally smooth; hence
$\IC(\Sym^d(\BA^2)_\tfP)=\uC[2\sK(\fP)]$, and
$\IC(\fM^\alpha_{\overline{\gamma,\emptyset,\fP}})=
N^{\gamma,\fP}_*\left(\IC(\fM^\gamma)\boxtimes\uC[2\sK(\fP)]\right)$.

Now we use Corollary ~\ref{decomp}, Proposition ~\ref{stalk}, and
induction in $\alpha,d$ (cf. the argument in ~\cite{b}, ~3.7).
\qed

\subsection{Uhlenbeck flag spaces for untwisted affine Lie algebras}
\label{speculation}
Let $\fg$ be a simple finite dimensional Lie algebra, and let
$\hg$ be the corresponding untwisted affine Lie algebra with the coroot
lattice $Y=\BZ[I]$, and the dual lattice of weights $X$. Let $\chg$
be the {\em Langlands dual} affine Lie algebra, with the roles of $X$ and $Y$
interchanged (note that if $\fg$ is not simply laced, then $\chg$ is twisted).
Let $\chn_+\subset\chg$ be the standard maximal nilpotent subalgebra graded
by $\BN[I]$, and $\chn(\fp)\subset\chn_+$ the nilpotent radical of the
standard maximal parabolic subalgebra $\fp\subset\chg$.
Let $\delta_0\in\BN[I]$ be the minimal imaginary root of $\chn_+$. Then
$\chn(\fp)$ is naturally graded by $\BN[\delta_0]:\ \chn(\fp)=\oplus_{d\in\BN}
\chn(\fp)_{d\delta_0}$. Also, $\chn(\fp)$ carries a natural integrable
action of the Langlands dual algebra $\fg^\vee$. Let $f\in\fg^\vee$ be a
principal nilpotent element. Let $W_\bullet\chn(\fp)$ be the monodromy
filtration associated to the action of $f$. Then the invariants
$(\chn(\fp))^f$ project injectively into $\gr_{W_\bullet}\chn(\fp)$, and
hence we obtain a grading on the space $(\chn(\fp))^f=:\fu(\chg)=
\oplus_{k\in\BN}\fu(\chg)^k$. Recall that we also have another grading
on $\fu(\chg)$, so that it is actually bigraded:
$\fu(\chg)=\bigoplus\limits_{d\in\BN}^{k\in\BN}\fu(\chg)^k_{d\delta_0}$.
Thus, the symmetric algebra $\Sym(\fu(\chg))$ is also bigraded:
$\Sym(\fu(\chg))=\bigoplus\limits_{d\in\BN}^{k\in\BN}
\Sym(\fu(\chg))^k_{d\delta_0}$.

\medskip

The Kashiwara definition of the flag scheme $\CB$, and the Drinfeld
definition of the based quasimaps' scheme $\bM^\alpha$ works for the
affine Lie algebra $\hg$ as well. Repeating the constructions
of ~\ref{family}, ~\ref{main char} we define the Uhlenbeck flag 
space $\fM^\alpha$. It is stratified as in ~\ref{uhl}. On the base
of ~\cite{bgfm}, ~7.3, we propose the following generalization of
Theorem ~\ref{IC} to arbitrary $\hg$.

\begin{conjecture}
\label{conj}
The stalk of $\IC(\fM^\alpha)$ at a point
$(\phi,\underline{c},\underline{s})\in\fM^\alpha_{\gamma,\Gamma,\fP}$
is isomorphic to
$$\bigotimes_{l=1}^m\left(\oplus_{r\in\BN}
\Sym^r(\chn_+)_{\beta_l}[2r]\right)^{\otimes
n_l}\otimes\bigotimes_{l=1}^g\left(\oplus_{r\in\BN}\Sym(\fu(\chg))^r_{d_l}
[2r]\right)^{\otimes k_l}[2|\gamma|]$$
\end{conjecture}

\section{Hecke correspondences}

\subsection{Boundary}
\label{bound}
We define an open subvariety $\CM^\alpha\supset\vM^\alpha\supset\oM^\alpha$
formed by the parabolic sheaves $\CF_\bullet$ which are locally free
parabolic sheaves in some Zariski open neighbourhood of $\bD_0\subset\CS'$.
The complementary closed subvariety $\CM^\alpha-\vM^\alpha$ is denoted by
$\nM^\alpha$. For $\alpha,\gamma\in\BN[I]$ we consider the {\em Hecke
correspondence} $\fE_\alpha^\gamma\subset\CM^\alpha\times\CM^{\alpha+\gamma}$
formed by the pairs $(\CF_\bullet,\CF'_\bullet)$ of parabolic sheaves such
that $\CF'_\bullet\subset\CF_\bullet$. The first projection 
$\fE_\alpha^\gamma\to\CM^\alpha$ will be denoted by $\bp$, and the second
projection $\fE_\alpha^\gamma\to\CM^{\alpha+\gamma}$ will be denoted by $\bq$.
Note that $\bq$ is proper while $\bp$ is not.

\begin{proposition}
\label{connect}
$\CM^\alpha$ is connected.
\end{proposition}

{\sl Proof:} Assume for a moment that $\vM^\alpha$ is connected for any
$\alpha$. Then we just have to prove that the boundary $\nM^\alpha$ is
connected. By induction in $\alpha$ we may assume that $\CM^\beta$ is
connected for any $\beta<\alpha$. Recall that $\alpha=\sum_{i\in I}a_ii$,
and for $i\in I$ such that $a_i>0$ we set $\alpha_i=\alpha-i$.
Then, evidently, $\nM^\alpha=\bigcup_{i}\bq(\fE_{\alpha_i}^i)$, and
any two nonempty pieces of the boundary intersect nontrivially:
$\bq(\fE_{\alpha_i}^i)\bigcap\bq(\fE_{\alpha_j}^j)\ne\emptyset$.
So we only have to prove that $\bq(\fE_{\alpha_i}^i)$ is connected.
But the fibers of projection $\fE_{\alpha_i}^i\to\CM^{\alpha_i}$ were
computed in ~\cite{fk2}, ~4.3.5. It follows in particular that these fibers
are connected. Since $\CM^{\alpha_i}$ is connected by induction assumption,
$\fE_{\alpha_i}^i$ is connected itself, hence $\bq(\fE_{\alpha_i}^i)$ is
also connected. Thus it remains to prove that $\vM^\alpha$ is connected.
This is the subject of the following Lemma.

\begin{lemma}
$\vM^\alpha$ is connected.
\end{lemma}

{\sl Proof:} We introduce a still bigger open subvariety
$\CM^\alpha\supset\cM^\alpha\supset\vM^\alpha$ formed by the parabolic
sheaves $\CF_\bullet$ such that $\CF_0$ is locally free in some Zariski
open neighbourhood of $\bD_0\subset\CS'$, and for $-n\leq k\leq0$
the quotient sheaves $\CF_k/\CF_{-n}$ on $\bD_0$ are locally free.
For such a parabolic sheaf $\CF_\bullet$ the quotients
$\CF_k/\CF_{-n}$ are locally free subsheaves of the vector bundle
$\CF_0|_{\bD_0}$ on $\bC$, and the parabolic sheaf $\CF_\bullet$ can be
uniquely reconstructed from $\CF_0$, and the flag of locally free
subsheaves $\CF_k/\CF_{-n}$ of $\CF_0|_{\bD_0}$.

Thus we have a cartesian diagram
$$
\begin{CD}
\cM^\alpha @>{\mathfrak r}>> \sQ^{\overline\alpha} \\
@V{r}VV            @V{\rho}VV       \\
\csM^{a_0} @>{\mathsf r}>> \Bun
\end{CD}
$$
Here $\csM^{a_0}\subset\CA^{a_0}$ is an open subvariety formed by
torsion free sheaves which are locally free in a Zariski open neighbourhood
of $\bD_0\subset\CS'$. Furthermore,
$\Bun$ is the moduli stack of $SL_n$-bundles on $\bC$, and
${\mathsf r}$ sends $\csM^{a_0}\ni\CF_0$ to $\CF_0|_{\bD_0}$.
Also, $r$ sends $\CF_\bullet$ to $\CF_0$.
Finally, $\sQ^{\overline\alpha}\stackrel{\rho}{\to}\Bun$ is Laumon's
stack of quasiflags of degree $\overline\alpha$, see ~\cite{la}.
Here $\overline\alpha$ is an element of coroot lattice $\BZ[I-0]\subset
\BZ[I]=Y$ of $\fsl_n\subset\hsl_n$. It is given by the formula
$\overline{\alpha}=\sum_{0\ne i\in I}(a_i-a_0)i$.

Now $\csM^{a_0}$ is connected being an open subvariety in $\CA^{a_0}$
which is connected by Nakajima's quiver realization ~\cite{n1}, ~2.1,
and cohomology computation ~\cite{nk}, ~2.3. Moreover, the fibers of
$\rho$ are connected since $\sQ^{\overline\alpha}$ is connected, and
$\Bun$ is normal. Hence $\cM^\alpha$ is connected. This completes the
proof of the Lemma along with Proposition ~\ref{connect}.
\qed

\begin{corollary}
\label{dimensio}
$\dim\CM^\alpha=2|\alpha|$.
\end{corollary}

{\sl Proof:} 
We use the cartesian diagram in the proof of the above Lemma, together with
the known formulas for the dimensions of $\Bun, \sQ^{\overline\alpha}$
(see ~\cite{la}) and of $\csM^{a_0}$ (see ~\cite{n1}).
\qed

\subsection{Generically trivial parabolic sheaves}
\label{pc_of_ma}

\nc{\br}{{\mathbf{r}}}
\nc{\sZ}{{\mathsf{Z}}}
\nc{\BO}{{\mathbb{O}}}
\nc{\BQ}{{\mathbb{Q}}}
\nc{\TCM}{{\widetilde{\CM}}}
%\nc{\bA}{{{\mathbf{A}}}}
\nc{\BF}{{{\mathbb{F}}}}
\nc{\fH}{{{\mathfrak{H}}}}
\nc{\tbH}{{\widetilde{\mathbf{H}}}}
\nc{\TA}{{\widetilde{\mathbf{A}}}}
\nc{\CMPC}[1]{{{\CM_{\text{\bf gt}}^{#1}}}}
\nc{\CMPCO}[1]{{{\overset{\circ}{\CM}_{\text{\bf gt}}^{#1}}}}
\nc{\fEF}{{{\mathfrak{E}\mathfrak{F}}}}
\nc{\fFE}{{{\mathfrak{F}\mathfrak{E}}}}
\nc{\fEO}{{{\overset{\circ}{\fE}{}}}}
\nc{\fEC}{{{\widehat{\fE}}{}}}
\nc{\ev}{{\mathsf{ev}}}
\nc{\lbb}{\{\!\{}
\nc{\rbb}{\}\!\}}
\nc{\tkappa}{{\tilde{\kappa}}}
%\nc{\tmu}{{\tilde{\mu}}}
\nc{\NR}{{\mathop{{\bf NR}}\nolimits}}
\renewcommand{\Im}{\mathop{{\sf Im}}\nolimits}

Our next goal is to study the action of Hecke correspondences $\fE^i$
on the cohomology of $\CM^\alpha$. Recall that the first projection 
$\bp:\ \fE_\alpha^i\to\CM^\alpha$ is not proper, and this causes a 
difficulty in the definition of the desired action. To get around this
difficulty we will introduce another version of moduli spaces 
$\CM^\alpha_{\mathbf{gt}}$
and Hecke correspondences between them which have proper projections.
Recall that $\CM^\alpha$ is the moduli space of parabolic sheaves which are
trivialized at $\bC\times\bx\bigcup\bc\times\bX$. In the definition of
$\CM^\alpha_{\mathbf{gt}}$ 
we replace the condition of triviality at $\bc\times\bX$
by the condition of triviality of $\CF_0$ at {\em some} line $c\times\bX$,
i.e. the condition of {\em generic triviality} of $\CF_0$. To give a
rigoruous definition we need some preparations. 

\medskip

For any collection of points $c_1,\dots,c_m$ of the curve $\bC$
we consider the moduli scheme $\TCM^\alpha(c_1,\dots,c_m)$ of all
infinite flags $\dots\subset\CF_{-1}\subset\CF_0\subset\CF_1\subset\dots$
of torsion free coherent sheaves of rank~$n$ on~$\CS'$ such that

(a) $\CF_{k+n}=\CF_k(\bD_0)$ for any $k$;

(b) $ch_1(\CF_k)=k[\bD_0]$ for any $k$: the first Chern classes are
proportional to the fundamental class of $\bD_0$;

(c) $ch_2(\CF_k)=a_i$ for $i\equiv k\pmod{n}$;

(d) For all $k\in\BZ$ the sheaf $\CF_k$ is locally free at
the lines $c_1\times\bX,\dots,c_m\times\bX\subset\CS'$;

(e) The sheaf $\CF_0$ is trivialized at the line $\bC\times\bx\subset\CS'$
and trivial at the lines $c_1\times\bX,\dots,c_m\times\bX\subset\CS'$.

\medskip

It is instructive to compare this definition with \ref{parab}.
The difference is the following. First, we replace one fixed
line $\bc\times\bX\subset\CS'$ with an $m$-tuple of lines
$c_1\times\bX,\dots,c_m\times\bX\subset\CS'$; and secondly,
we drop the condition (e) of \ref{parab}, imposing the behavior
of the restrictions ${\CF_k}_{|\bc\times\bX}$ and replaced it
by a much weaker condition of locally freeness.

\begin{proposition}
\label{tcma}
For all $m>0$ and all $c_1,\dots,c_m\in C$ the moduli scheme
$\TCM^\alpha(c_1,\dots,c_m)$ exists. It is smooth connected
scheme of dimension $\dim\sB+|2\alpha|$. For any permutation
$\sigma$ of the set $\{1,\dots,m\}$ the schemes
$\TCM^\alpha(c_1,\dots,c_m)$ and
$\TCM^\alpha(c_{\sigma(1)},\dots,c_{\sigma(m)})$
are canonically isomorphic.
\end{proposition}

{\sl Proof:}
The first part follows from~\cite{hl} and~\cite{y} as in~\ref{CM},
and the third part is evident.
Thus it remains to check the smoothness and the connectedness
and to compute the dimension.

First, consider the moduli scheme $\TCM^\alpha(\bc)$.
Consider the evaluation map
$$
\ev_{(\bc,\by)}:\TCM^\alpha(\bc) \to \sB,\quad
(\CF_k)\mapsto
\big(\Im({\CF_k}_{|(\bc,\by)}\to{\CF_0}_{|(\bc,\by)})\big)_{-n\le k\le 0}.
$$
The locally freeness condition \ref{pc_of_ma} (d) implies that
${\CF_k}_{|(\bc,\by)}$ is an $n$-dimensional vector space, while
the conditions (a) and (b) imply that the map
${\CF_k}_{|(\bc,\by)}\to{\CF_0}_{|(\bc,\by)}$ has rank $k+n$.
Hence $\ev_{(\bc,\by)}(\CF_k)$ is a flag in the vector space
${\CF_0}_{(\bc,\by)}$ which is canonically isomorphic to $V$,
thus a point of the flag variety $\sB$.
Now note that the map $\ev_{(\bc,\by)}$ is evidently a locally
trivial fibration, and the fiber of $\ev_{(\bc,\by)}$ is canonically
isomorphic to the variety $\CM^\alpha$. Hence Lemma~\ref{yoko} and
Remark~\ref{dimcm} imply that $\TCM^\alpha(\bc)$ is a smooth connected
variety of dimension $\dim\sB+2|\alpha|$.

Further, for any point $c_1\in\bC$ choose an automorphism $g$ of
the curve $\bC$ such that $g(c_1)=\bc$. Then $g$ identifies the
moduli schemes $\TCM^\alpha(c_1)$ and $\TCM^\alpha(\bc)$, hence
$\TCM^\alpha(c_1)$ is a smooth connected variety of dimension
$\dim\sB+2|\alpha|$ for all $c_1\in\bC$.

Finally, it is clear that $\TCM^\alpha(c_1,\dots,c_m)$ is an open
subscheme of $\TCM^\alpha(c_1)$, hence $\TCM^\alpha(c_1,\dots,c_m)$
is a smooth connected variety of dimension $\dim\sB+2|\alpha|$ for all
$c_1,\dots,c_m\in\bC$.
\qed

\begin{definition}
Let $\CMPC\alpha$ denote the gluing of schemes $\TCM^\alpha(c_1)$
for all $c_1\in\bC$ with respect to the open subsets
$\TCM^\alpha(c_1)\supset\TCM^\alpha(c_1,c_2)\subset\TCM^\alpha(c_2)$.
\end{definition}

\begin{theorem}
\label{cmpc}
The scheme $\CMPC\alpha$ is a smooth connected variety of
dimension $\dim\sB + 2|\alpha|$. The moduli schemes $\TCM^\alpha(c_1)$
form an open covering of $\CMPC\alpha$ and
$\TCM^\alpha(c_1)\cap\TCM^\alpha(c_2) = \TCM^\alpha(c_1,c_2)$.
\end{theorem}

{\sl Proof:}
The only thing we need to check is that $\CMPC\alpha$ is
a scheme of finite type. Then all the rest follows from
Proposition~\ref{tcma}. On the other hand, if we want to
check that $\CMPC\alpha$ is of finite type, it suffices
to show that there exists an integer $s$ such that for any
collection of distinct points $c_1,c_2^1,\dots,c_2^s\in\bC$
we have
$$
\TCM^\alpha(c_1) = \bigcup_{p=1}^s\TCM^\alpha(c_1,c_2^p).
$$
Then it would follow that $\CMPC\alpha$ is in fact a gluing
of $(s+1)$ moduli schemes $\TCM^\alpha(c_1^i)$ for an arbitrary
collection of $(s+1)$ distinct points $c_1^1,\dots,c_1^{s+1}\in\bC$.

Let us show that $s=|\alpha|+1$ is big enough. Since the group
of automorphisms of $\bC$ acts transitively, it suffices to
consider only the case $c_1=\bc$. So let $\CF_k$ be a point of
$\TCM^\alpha(\bc)$. Let $\CN(\CF)_k$ be its saturation and
denote $\beta=\sum b_i i$, where $b_k=ch_2(\CN(\CF)_k)$.
Then $\CN(\CF)_0\in\CA^{a_0-b_0}$ and for any $k\in\BZ$
the sheaf $\CN(\CF)_k/\CF_k$ is a sheaf on $\CS'$ of length $b_k$.
Recall that according to Lemma ~\ref{B12} the subset $D\subset(\bC-\bc)$
formed by the points $x\in\bC$ such that the restriction 
$\CN(\CF)_0|{x\times\bX}$ is nontrivial, has cardinality at most $a_0-b_0$.
Therefore
$$
D\cup
\left(\bigcup_{k=0}^{n-1}\varrho_\bC(\supp(\CN(\CF)_k/\CF_k))\right)
$$
is a subset in $\bC-\bc$ of cardinality not greater than
$$
(a_0-b_0) + \sum_{k=0}^{n-1} b_k = (a_0-b_0) + |\beta|.
$$
Now if $s>(a_0-b_0) + |\beta|$ and $c_2^1,\dots,c_2^s$ are distinct
points of $\bC-\bc$ then there exists $1\le p\le s$ such that for all
$k\in\BZ$ the sheaf $\CF_k$ is locally free and the sheaf $\CF_0$
is trivial at the line $c_2^p\times\bX$. Thus $\CF_\bullet$ lies in
$\TCM^\alpha(\bc,c_2^p)$. Thus it remains to check that the integer
$(a_0-b_0) + |\beta|$ is uniformly bounded for all $\CF_\bullet$.
But this is evident, because we always have $\beta\le\alpha$, hence
$(a_0-b_0) + |\beta| \le |\alpha$. Thus $s=|\alpha|+1$ is indeed
big enough.
\qed

%\subsection{Cohomology of $\CMPC\alpha$}
%\label{coh_cmpc}

%Let $P(V^\bullet,t)$ denote the Poincare polynomial of
%a graded vector space $V^\bullet$. Then it is easy to show that
%$$
%P(H^\bullet(\CMPC\alpha,\BQ),t) =
%\sum_{\emptyset\ne J\subset \{1,\dots,s\}}
%(-1)^{|J|}P(H^\bullet(\TCM^\alpha(c_J),\BQ),t),
%$$
%where $s=|\alpha|+1$, $c_1,\dots,c_s\in C$ is a collection
%of pairwise distinct points, and $c_J$ denotes the subcollection
%given by an index set $J$. It is easy to show that
%$$
%P(H^\bullet(\TCM^\alpha(c_1),\BQ),t) =
%t^{|\alpha|}\left(\sum_{\kappa\in\fK_{\hgl_n}}t^{-\sK(\kappa)}\right)
%\left(\sum_{w\in W}t^{\ell(w)}\right),
%$$
%where $W$ is the Weyl group of $SL_n$ and $\ell(w)$ is the length function.
%However, I don't know how to compute $P(H^\bullet(\TCM^\alpha(c_J),\BQ),t)$
%in the case when $|J|>1$.

\subsection{Correspondences}
\label{cors}

For any $\alpha,\gamma\in\BN[I]$ we define the {\em Hecke correspondence}
$\fE_\alpha^\gamma\subset\CMPC\alpha\times\CMPC{\alpha+\gamma}$ as a closed
subvariety formed by all the pairs $(\CF_\bullet,\CF'_\bullet)\in
\CMPC\alpha\times\CMPC{\alpha+\gamma}$ such that
$\CF'_\bullet\subset\CF_\bullet$, and the quotient is supported at
$\bD_0\subset\CS'$.

We have canonical projections
$$
\bp:\fE_\alpha^\gamma\to\CMPC\alpha,\quad
\bq:\fE_\alpha^\gamma\to\CMPC{\alpha+\gamma},\quad\text{and}\quad
\br:\fE_\alpha^\gamma\to\bC^\gamma.
$$
Here $\bp$ and $\bq$ are induced by projections of the product
$\CMPC\alpha\times\CMPC{\alpha+\gamma}$ to the factors,
and $\br$ is defined as follows:
$$
\br(\CF_\bullet,\CF'_\bullet) =
\supp(\CF_\bullet,\CF'_\bullet) =
\sum_i \supp(\CF_i/\CF'_i)i \in \bD_0^\gamma = \bC^\gamma.
$$

\begin{lemma}
\label{proper}
The maps $\bp$ and $\bq$ are proper.
\end{lemma}
{\sl Proof:}
Evident.
\qed

\subsection{Top-dimensional components}

We begin with some notation. Recall that for any
$\CF_\bullet\subset\CF'_\bullet$ the quotient
$T_\bullet = \CF_\bullet/\CF'_\bullet$ can be considered
as a representation of the infinite linear quiver $\bA_\infty$
in the category of torsion sheaves supported on $\bD_0$.
It is clear that $T_\bullet$ is a nilpotent representation.
On the other hand, the periodicity of $\CF_\bullet$ and $\CF'_\bullet$
imply the periodicity of $T_\bullet$, namely a canonical isomorphism
(the triviality of the normal bundle $\CN_{\bD_0/\CS'}$ is used here)
$$
T_{k+n} \cong T_k.
$$
Thus $T_\bullet$ can (and will) be considered as a nilpotent
representation of the cyclic quiver $\TA_{n-1}$. Following~\cite{fk2}
we denote by $\NR_n(\bD_0)$ the category of nilpotent representations
of the cyclic quiver $\TA_{n-1}$ in the category of torsion sheaves
supported on $\bD_0$. For any object $T_\bullet\in\NR_n(\bD_0)$
and a point $x\in\bD_0$ we denote by $\Gamma_x(T_\bullet)$
the representation of the cyclic quiver $\TA_{n-1}$ in the
category of vector spaces, formed by sections of $T_\bullet$
with support at the point $x\in\bD_0$. Recall that the isomorphism
classes of nilpotent representations of the cyclic quiver $\TA_{n-1}$
in the category of vector spaces are numbered by Kostant partitions
of $\hgl_n$. We denote by $\kappa_x(T_\bullet)\in\fK_{\hgl_n}$
the isomorphism class of $\Gamma_x(T_\bullet)$.

Now we are going to use the results of \cite{fk2} to describe
the top-dimensional components of $\fE_\alpha^\gamma$. Choose a Kostant
partition $\kappa=\lbb\theta_1,\dots,\theta_m\rbb\in\fK_{\hgl_n}(\gamma)$,
where $\theta_1,\dots,\theta_m\in R^+_{\hgl_n}$. Consider a subset
$\fEO{}_\alpha^\kappa \subset \fE_\alpha^\gamma$ consisting of all
pairs $(\CF_\bullet,\CF'_\bullet)$ such that

\begin{enumerate}
\item $\br(\CF_\bullet,\CF'_\bullet)\in \bC^\gamma_\kappa$, that is
$\supp(\CF_\bullet,\CF'_\bullet)=\sum_{r=1}^m|\theta_r|x_r$ with
all $x_r$ distinct;

\item $\CF_\bullet$ is locally free at the points $x_1,\dots,x_m$;

\item $\kappa_{x_r}(\CF_\bullet/\CF'_\bullet)=\lbb\theta_r\rbb$
for all $1\le r\le m$.
\end{enumerate}

We define $\fE_\alpha^\kappa$ as the closure of $\fEO{}^\kappa_\alpha$.

\begin{proposition}
\label{ic_fe}
Dimension of any irreducible component of $\fE_\alpha^\gamma$
is not greater than $\dim\sB + 2|\alpha| + |\gamma|$. Any component
of this dimension coincides with $\fE_\alpha^\kappa$ for
some $\kappa\in\fK_{\hgl_n}(\gamma)$.
\end{proposition}

{\sl Proof:}
Consider a stratification of $\CMPC\alpha\times \bC^\gamma$
via the defect of $\CF_\bullet$ at the support of
$\sum\gamma_rx_r\in \bC^\gamma$, namely
$$
\CMPC\alpha\times \bC^\gamma =
\bigsqcup_{\begin{Sb}
\Gamma\in\Gamma(\gamma)\\
|\kappa'_1|+\dots+|\kappa'_m|=\gamma'\le\alpha
\end{Sb}}
\sZ^\Gamma_\alpha(\kappa'_1,\dots,\kappa'_m).
$$
Here $\sZ^\Gamma_\alpha(\kappa'_1,\dots,\kappa'_m)\subset
\CMPC\alpha\times \bC^\gamma$ is the subspace of all pairs
$(\CF_\bullet,\sum\gamma_rx_r)$ such that
\begin{enumerate}
\item $\{\{\gamma_1,\dots,\gamma_m\}\}=\Gamma$;
\item $\kappa_{x_r}(\CN(\CF)_\bullet/\CF_\bullet)=\kappa'_r$
for all $1\le r\le m$;
\end{enumerate}

Consider the partial saturation map
$$
\sZ^\Gamma_\alpha(\kappa'_1,\dots,\kappa'_m)\to
\CMPC{\alpha-\gamma'}\times \bC^\gamma_\Gamma,\quad
(\CF_\bullet,\sum\gamma_rx_r)\mapsto
(\CN_{x_1,\dots,x_m}(\CF)_\bullet,\sum\gamma_rx_r).
$$
The fiber of this saturation map was described in \cite{fk2}. It was
denoted there by $K_\mu$, where $\mu=\{\{\kappa'_1,\dots,\kappa'_m\}\}$ ---
the corresponding multipartition. It was shown in Lemma~3.1.4 and Theorem~1
of {\em loc.\ cit.} that
$$
\dim K_\mu = \sum(||\kappa'_r||-\sK(\kappa'_r)).
$$
This implies
$$
\arraycolsep=2pt
\begin{array}{rcl}
\dim \sZ^\Gamma_\alpha(\kappa'_1,\dots,\kappa'_m) &=&
\dim\sB + 2|\alpha-\gamma'| + m + \sum(||\kappa'_r||-\sK(\kappa'_r))
\smallskip\\ &=&
\dim\sB + 2|\alpha|-|\gamma'| - \sum\sK(\kappa'_r) + m
\smallskip\\ &=&
\dim\sB + 2|\alpha| + |\gamma| + \sum(1-|\gamma_r+\gamma'_r|-\sK(\kappa'_r)),
\end{array}
$$
where $\gamma'_r=|\kappa'_r|$.

Now consider a stratification of $\fE_\alpha^\gamma$
via the defect of the sheaves $\CF_\bullet$ and $\CF'_\bullet$
at the support of $\CF'_\bullet/\CF_\bullet$, namely
$$
\fE_\alpha^\gamma =
\bigsqcup_{\begin{Sb}
\Gamma\in\Gamma(\gamma)\\
|\kappa'_1|+\dots+|\kappa'_m|=\gamma'\le\alpha\\
|\tkappa_r|=|\kappa'_r|+\gamma_r
\end{Sb}}
\sZ^\Gamma_\alpha(\kappa'_1,\dots,\kappa'_m;\tkappa_1,\dots,\tkappa_m).
$$
Here $\sZ^\Gamma_\alpha(\kappa'_1,\dots,\kappa'_m;\tkappa_1,\dots,\tkappa_m)
\subset\fE_\alpha^\gamma$ is the subspace of all pairs
$(\CF_\bullet,\CF'_\bullet)$ such that
\begin{enumerate}
\item $\br(\CF_\bullet,\CF'_\bullet)=\sum\gamma_rx_r\in \bC^\gamma_\Gamma$;
\item $\kappa_{x_r}(\CN(\CF)_\bullet/\CF_\bullet)=\kappa'_r$
for all $1\le r\le m$;
\item $\kappa_{x_r}(\CN(\CF')_\bullet/\CF'_\bullet)=\tkappa_r$
for all $1\le r\le m$;
\end{enumerate}
Note that $\CF'_\bullet\subset\CF_\bullet$ implies $\CN(\CF')=\CN(\CF)$,
hence we indeed have $|\tkappa_r|=|\kappa'_r|+\gamma_r$.

Now consider the map $\bp\times\br$ restricted to the stratum
$\sZ^\Gamma_\alpha(\kappa'_1,\dots,\kappa'_m;\tkappa_1,\dots,\tkappa_m)$.
It is clear that
$$
(\bp\times\br)
(\sZ^\Gamma_\alpha(\kappa'_1,\dots,\kappa'_m;\tkappa_1,\dots,\tkappa_m))
= \sZ^\Gamma_\alpha(\kappa'_1,\dots,\kappa'_m).
$$
Furthermore, it is easy to see that the fiber of this map
over a point $(\CF_\bullet,\sum\gamma_rx_r)$ consists
of all $\CF'_\bullet$, such that
\begin{enumerate}
\item $\kappa_{x_r}(\CN(\CF)_\bullet/\CF'_\bullet)=\tkappa_r$
for all $1\le r\le m$;
\item $\CF'_\bullet\subset\CF$.
\end{enumerate}
It follows that this fiber can be embedded into
the variety $K_{\tmu}$, where $\tmu=\lbb\tkappa_1,\dots,\tkappa_m\rbb$
as a closed subvariety (the closed condition is the condition 2 above).
In particular, the dimension of the fiber is not greater than
$$
\sum_{r=1}^m(||\tkappa_r||-\sK(\tkappa_r)) =
\sum_{r=1}^m(|\gamma'_r+\gamma_r| - \sK(\tkappa_r)).
$$
Comparing this with the formula for the dimension of the
stratum $\sZ^\Gamma_\alpha(\kappa'_1,\dots,\kappa'_m)$ we see that
\begin{multline*}
\dim\sZ^\Gamma_\alpha(\kappa'_1,\dots,\kappa'_m;\tkappa_1,\dots,\tkappa_m)
\\ \le
\dim\sB + 2|\alpha| + |\gamma| +
\sum(1-|\gamma_r+\gamma'_r|-\sK(\kappa'_r)) +
\sum_{r=1}^m(|\gamma'_r+\gamma_r| - \sK(\tkappa_r))
\\ =
\dim\sB + 2|\alpha| + |\gamma| +
\sum_{r=1}^m(1-\sK(\kappa'_r)-\sK(\tkappa_r)).
\end{multline*}
Note that since $\gamma_r>0$ we have $\sK(\tkappa_r)\ge 1$,
hence the last term is always nonpositive. Therefore,
$$
\dim\sZ^\Gamma_\alpha(\kappa'_1,\dots,\kappa'_m;\tkappa_1,\dots,\tkappa_m) \le
\dim\sB + 2|\alpha| + |\gamma|
$$
and the equality is possible only when $\sK(\kappa'_r)=0$, $\sK(\tkappa_r)=1$
for all $1\le r\le m$. This means that $\CF_\bullet$ is locally free at
the points $x_1,\dots,x_r$ and that $\tkappa_r=\lbb\theta_r\rbb$ for some
$\theta_r\in R^+_{\hgl_n}$. Moreover, it is easy to see that
in the latter case the condition (2) above is void, hence
$$
\dim\sZ^\Gamma_\alpha(0,\dots,0;\lbb\theta_1\rbb,\dots,\lbb\theta_1\rbb) =
\dim\sB + 2|\alpha| + |\gamma|.
$$
Finally, it remains to note that
$\fEO{}_\alpha^{\lbb\theta_1,\dots,\theta_m\rbb} =
\sZ^\Gamma_\alpha(0,\dots,0;\lbb\theta_1\rbb,\dots,\lbb\theta_1\rbb)$.
\qed

\subsection{}

In Proposition~\ref{ic_fe} we gave an explicit description
of open parts of the top-dimensional irreducible components
of $\fE_\alpha^\gamma$. Below we will need for technical reasons
also an explicit description of some closed subset
$\fEC{}_\alpha^\kappa\subset\fE_\alpha^\gamma$ such that
$\fE_\alpha^\kappa\subset\fEC{}_\alpha^{\kappa'}$ iff $\kappa=\kappa'$.

Now we will define such closed subsets. We begin with some notation.
For $T_\bullet\in\NR_n(\bD_0)$ we denote
$$
\Gamma(T_\bullet)=\Gamma(\CS,T_\bullet)=\oplus_{x\in\bD_0}\Gamma_x(T_\bullet),
\quad\text{and}\quad
\kappa(T_\bullet)=\sum_{x\in\bD_0}\kappa_x(T_\bullet).
$$
Thus $\kappa(T_\bullet)$ is the isomorphism class of $\Gamma(T_\bullet)$.
Now, for every $\gamma\in\BZ[I]$ let $\bV(\gamma)$ denote the
representation space of all $\gamma$-dimensional nilpotent
representations of the cyclic quiver $\TA_{n-1}$ and let $\GL(\gamma)$
denote the group acting on $\bV(\gamma)$ by change of bases.
Then the $\GL(\gamma)$-orbits on $\bV(\gamma)$ are nothing but
the isomorphism classes of nilpotent representations of $\TA_{n-1}$.
Let $\BO_\kappa$ denote the orbit corresponding to a Kostant partition
$\kappa$. Note that we have a canonical partial order on the set of
orbits. It induces a partial order on the set of Kostant partition:
for $\kappa,\kappa'\in\fK_{\hgl_n}(\gamma)$ we have
$$
\kappa \le \kappa'\quad\text{iff}\quad
\BO_\kappa\subset\overline{\BO_{\kappa'}}\subset\bV(\gamma).
$$

Now take an arbitrary Kostant partition $\kappa\in\fK_{\hgl_n}(\gamma)$
and consider the subset $\fEC{}_\alpha^\kappa\subset\fE_\alpha^\gamma$
defined as
$$
\fEC{}_\alpha^\kappa = \{(\CF_\bullet\supset\CF'_\bullet)\ |\
\kappa(\CF_\bullet/\CF'_\bullet)\le \kappa\text{ and }
\br(\CF_\bullet,\CF'_\bullet)\in \overline{\bC^\gamma_\kappa} \}
$$

\begin{lemma}
\label{fec}
The subset $\fEC{}_\alpha^\kappa\subset\fE_\alpha^\gamma$ is closed and
$\fEO{}_\alpha^\kappa\subset\fEC{}_\alpha^{\kappa'}$ if and only if
$\kappa=\kappa'$.
\end{lemma}

{\sl Proof:} It is clear that both conditions defining
$\fEC{}_\alpha^\kappa\subset\fE_\alpha^\gamma$ are closed, hence
$\fEC{}_\alpha^\kappa$ is closed. Now assume that
$\fEO{}_\alpha^\kappa\subset\fEC{}_\alpha^{\kappa'}$.
It follows than that
$$
\kappa \le \kappa'\qquad\text{and}\qquad
\bC^\gamma_\kappa\subset\overline{\bC^\gamma_{\kappa'}}.
$$
Let us check that this is possible only if $\kappa=\kappa'$.
%In fact we will check a stronger assertion, we will show that
%the orders on the set of all Kostant partitions, given by adjunction
%of orbits $\BO_\kappa$ and by adjunction of strata $\bC^\gamma_\kappa$
%are opposite. To this end we use the following result of Guo~\cite{GUO}:
%the partial order $\le$ given by adjunction of orbits $\BO_\kappa$
%is compatible with the linear order defined as follows:
%$$
%\kappa < \kappa' \quad \text{\bf \dots to be inserted !}
%$$
It is clear that the partial order given by adjacency
of strata $\bC^\gamma_\kappa$ can be described as follows. Assume that
$\theta_1,\dots,\theta_r,\theta'_1,\dots,\theta'_m\in R^+_{\hgl_n}$,
$m=s_1+\dots+s_r$, and for all $1\le p\le r$ we have
$$
|\theta_p|=\sum_{q=1}^{s_p}|\theta'_{s_1+\dots+s_{p-1}+q}|.
$$
Then we have
$$
\bC^\gamma_{\lbb\theta_1,\dots,\theta_r\rbb} \subset
\overline{\bC^\gamma_{\lbb\theta'_1,\dots,\theta'_m\rbb}}
$$
and all adjacencies have such form. 
In particular, if $\overline{\bC^\gamma_\kappa}
\subset\overline{\bC^\gamma_{\kappa'}}$
is a proper inclusion, then the number of entries $\sK(\kappa)$ is
strictly smaller than $\sK(\kappa')$. 

On the other hand, Ringel's explicit description of the order $\kappa\leq
\kappa'$ in ~\cite{r}, ~4.7 implies that we must have $\sK(\kappa)\geq
\sK(\kappa')$ whenever $\kappa\leq\kappa'$.  
The Lemma follows.
\qed

\subsection{Action of the Hall algebra}
\label{hall}

Let $\tbH_n$ denote the generic Hall algebra of the category of
nilpotent representation of the cyclic quiver $\TA_{n-1}$.
Recall that the generic algebra $\tbH_n$ is an algebra over $\BZ[q]$
(polynomials in a formal variable $q$) with a basis $S_\kappa$
indexed by isomorphism classes of representations, that is by
Kostant partitions, and with the following multiplication rule
$$
S_{\kappa'}\cdot S_{\kappa''} =
\sum_\kappa c^\kappa_{\kappa',\kappa''}(q)S_\kappa,
$$
where the structure constants $c^\kappa_{\kappa',\kappa''}(q)$ are defined
as follows. Take $\BF_q$ for a base field and choose a representation
$W_\bullet$ over $\BF_q$ in the isomorphism class $\kappa$. Then
$c^\kappa_{\kappa',\kappa''}(q)$ is defined as the number of
subrepresentations $W'_\bullet\subset W_\bullet$, such that
the isomorphism class of $W'_\bullet$ equals $\kappa'$ and
the isomorphism class of $W_\bullet/W'_\bullet$ equals $\kappa''$.
It turns out that $c^\kappa_{\kappa',\kappa''}(q)$ is a polynomial
function of $q$, thus we can consider $\tbH_n$ as an algebra over $\BZ[q]$.

 From now on we consider a specialization of the Hall algebra $\tbH_n$
at $q=1$ and denote it by $\bH$.. As before, $\bH$ is a $\BQ$-algebra
having $S_\kappa$ for a basis and $c^\kappa_{\kappa',\kappa''}(1)$
for structure constants.

Consider a graded vector spaces
\begin{equation}
\label{fh}
\fH = \bigoplus_{\alpha\in\BZ[I]} H^{\bullet-|\alpha|}(\CMPC\alpha,\BQ)
%\quad\text{and}\quad
%\fH_c = \bigoplus_{\alpha\in\BZ[I]} H_c^{\bullet-|\alpha|}(\CMPC\alpha,\BQ)
\end{equation}
(note the shift of the cohomological grading).

For each Kostant partition $\kappa\in\fK_{\hgl_n}(\gamma)$ we consider
an operator on cohomology given by a correspondence
$\fE_\alpha^\kappa\subset\CMPC\alpha\times\CMPC{\alpha+\gamma}$:
$$
e_\kappa = [\fE_\alpha^\kappa]:
H^\bullet(\CMPC\alpha,\BQ) \to H^\bullet(\CMPC{\alpha+\gamma},\BQ).
$$
Since $\dim\CMPC{\alpha+\gamma}=\dim\CMPC\alpha + 2|\gamma|$ and
$\dim\fE_\alpha^\kappa=\dim\CMPC\alpha+|\gamma|$ it follows that
$e_\kappa$ shifts the cohomological degree by $|\gamma|$. Hence
$e_\kappa$ considered as an operator in the vector space $\fH$
preserves the cohomological degree.

\begin{remark}
\label{rem}
Instead of defining $e_\kappa$ as the operator given by the correspondence
$[\fE_\alpha^\kappa]$ we could define $e_\kappa$ as the component of
the operator given by the correspondence $[\fEC{}_\alpha^\kappa]$,
which increases the cohomological dimension by~$2||\kappa||$. According to
Lemma~\ref{fec} these definitions are equivalent.
\end{remark}

\begin{theorem}
\label{hallth}
The vector space $\fH$ is a module over the Hall algebra $\bH$,
where $S_\kappa$ act via $e_\kappa$.
\end{theorem}

{\sl Proof:}
Let $\kappa'\in\fK_{\hgl_n}(\gamma')$, $\kappa''\in\fK_{\hgl_n}(\gamma'')$
and put $\gamma=\gamma'+\gamma''$. We have to compute
the composition of correspondences
$[\fE_{\alpha+\gamma'}^{\kappa''}]\cdot[\fE_\alpha^{\kappa'}]$.
Instead, using Remark~\ref{rem} we can compute the component
of the composition of correspondences
$[\fEC{}_{\alpha+\gamma'}^{\kappa''}]\cdot[\fEC{}_\alpha^{\kappa'}]$
that increase the cohomological dimension by $2|\gamma|$.

Consider the product
$\CMPC\alpha\times\CMPC{\alpha+\gamma'}\times\CMPC{\alpha+\gamma}$
and let $p_{ij}$ denote the projection to the product of the $i$-th
and $j$-th factors. Consider the subset
$\fEC{}_\alpha^{\kappa',\kappa''}\subset
\CMPC\alpha\times\CMPC{\alpha+\gamma'}\times\CMPC{\alpha+\gamma}$
defined as
\begin{multline*}
\fEC{}_\alpha^{\kappa',\kappa''} :=
p_{12}^{-1}(\fEC{}_{\alpha}^{\kappa'}) \cap
p_{23}^{-1}(\fEC{}_{\alpha+\gamma'}^{\kappa''}) \\ =
\{(\CF_\bullet\supset\CF'_\bullet\supset\CF''_\bullet)\ |\
(\CF_\bullet,\CF'_\bullet)\in\fEC{}_{\alpha}^{\kappa'}\text{ and }
(\CF'_\bullet,\CF''_\bullet)\in\fEC{}_{\alpha+\gamma'}^{\kappa''} \}
\end{multline*}
Then
$[\fEC{}_{\alpha}^{\kappa'}]\cdot[\fEC{}_{\alpha+\gamma'}^{\kappa''}]$
is given by ${p_{13}}_*[\fEC{}_\alpha^{\kappa',\kappa''}]$.
But it is clear that
$p_{13}(\fEC{}_\alpha^{\kappa',\kappa''}) \subset \fE_\alpha^\gamma$,
hence by Proposition~\ref{ic_fe} the component of
$[\fEC{}_{\alpha+\gamma'}^{\kappa''}]\cdot[\fEC{}_\alpha^{\kappa'}]$
increasing the cohomological dimension by~$2|\gamma|$ equals to
$$
\sum_{\kappa\in\fK_{\hgl_n}} d^\kappa_{\kappa',\kappa''}
[\fE_\alpha^\kappa]
$$
for some constants $d^\kappa_{\kappa',\kappa''}$ which we have to compute.
Further, it is clear that $d^\kappa_{\kappa',\kappa''}$ equals the
number of points of $\fEC{}_\alpha^{\kappa',\kappa''}$
over a generic point of $\fE_\alpha^\kappa$. Since we are interested
in a generic point, we can take a point in $\fEO{}_\alpha^\kappa$.
So let $(\CF_\bullet,\CF''_\bullet)\in\fEO{}_\alpha^\kappa$ and
denote $T_\bullet=\CF_\bullet/\CF''_\bullet$. Then it is clear that
$d^\kappa_{\kappa',\kappa''}$ equals the number of subobjects
$T''_\bullet\subset T_\bullet$ such that for
$T'_\bullet=T_\bullet/T''_\bullet$ the following conditions are satisfied:

\begin{enumerate}
\item $\kappa(T'_\bullet)\le \kappa'$;
\item $\supp(T'_\bullet)\subset\overline{\bC^{\gamma'}_{\kappa'}}$;
\item $\kappa(T''_\bullet)\le \kappa''$;
\item $\supp(T''_\bullet)\subset\overline{\bC^{\gamma''}_{\kappa''}}$;
\end{enumerate}

Now assume that $\kappa=\lbb\theta_1,\dots,\theta_m\rbb$
and that $\supp(T_\bullet)=\sum\theta_rx_r$. Then it is
clear that $\kappa_{x_r}(T_\bullet)=\lbb\theta_r\rbb$ for all $1\le r\le m$.
Assume that $T''_\bullet$ is a subobject in $T_\bullet$.
Then $W''_r=\Gamma_{x_r}(T''_\bullet)\subset\Gamma_{x_r}(T_\bullet)=W_r$
is a subrepresentation. Moreover, $T''_\bullet$ is uniquely determined
by this collection of subrepresentations. Indeed, it is equal to
the image of the natural map
$$
\oplus_{r=1}^mW''_r\otimes\CO_\CS \to
\oplus_{r=1}^mW_r\otimes\CO_{\CS} \to
T_\bullet.
$$
Finally, note that the set of all nontrivial subrepresentations
$W''_r\subset W_r$ is in a bijection with the set of all $\theta''_r$
such that $\theta''_r$ ends at the same vertex as $\theta_r$ does,
and has smaller length. Put $\theta'_r=\theta_r/\theta''_r$. Then we have
$\theta_r=\theta'_r\star\theta''_r$. Note that if $T''_\bullet$ and
$T'_\bullet$ is the subobject and the quotient object of $T_\bullet$
corresponding to such collection $\theta''_r$ then
$$
\begin{array}{rclrcl}
\kappa(T'_\bullet) &=& \lbb\theta'_1,\dots,\theta'_m\rbb,\qquad&
\supp(T'_\bullet) &\in& C_{\gamma'}^{\lbb\theta'_1,\dots,\theta'_m\rbb};\\
\kappa(T''_\bullet) &=& \lbb\theta''_1,\dots,\theta''_m\rbb, &
\supp(T''_\bullet) &\in& C_{\gamma''}^{\lbb\theta''_1,\dots,\theta''_m\rbb};
\end{array}
$$
Thus $d^\kappa_{\kappa',\kappa''}$ equals the number of collections
$(\theta'_r,\theta''_r)_{r=1}^m$ such that
$$
\kappa'=\lbb\theta'_1,\dots,\theta'_m\rbb,\quad
\kappa''=\lbb\theta'_1,\dots,\theta'_m\rbb,\quad\text{and}\quad
\theta_r=\theta'_r\star\theta''_r \text{ for all } 1\le r\le m.
$$
It remains to note that this is precisely $c^\kappa_{\kappa',\kappa''}(1)$
(see e.g. ~\cite{fm}).
\qed

\subsection{Action of $\hsl_n$}

In addition to the operators $e_\kappa$ introduced above, we define
operators $f_\kappa$ as the operators on the cohomology induced
by the transposed correspondences:
$$
f_\kappa = [(\fE_{\alpha-\gamma}^\kappa)^T]:
H^\bullet(\CMPC\alpha,\BQ) \to H^\bullet(\CMPC{\alpha-\gamma},\BQ).
$$

To unburden the notation denote the operators $e_{\lbb i\rbb}$ and
$f_{\lbb i\rbb}$ by $e_i$ and $f_i$ respectively. Further, define
the operator $h_i$ on $H^\bullet(\CMPC\alpha,\BQ)$ as a scalar
$\langle i',\alpha\rangle+2$-multiplication.

\begin{proposition}
\label{ef}
We have $[e_i,f_j]=\delta_{ij}h_i$.
\end{proposition}

{\sl Proof:}
We have to compare the following compositions of correspondences:
$$
e_if_j = [\fE_{\alpha}^{\lbb i\rbb}]\cdot[\fE_{\alpha+i-j}^{\lbb j\rbb}]^T
\quad\text{and}\quad
f_je_i = [\fE_{\alpha-j}^{\lbb j\rbb}]^T\cdot[\fE_{\alpha-j}^{\lbb i\rbb}].
$$
Instead, as in the Proof of Theorem~\ref{hallth} we will compare
the components of the compositions
$$
[\fE_{\alpha}^i]\cdot[\fE_{\alpha+i-j}^j]^T\quad\text{and}\quad
[\fE_{\alpha-j}^j]^T\cdot[\fE_{\alpha-j}^i]
$$
preserving the cohomological degree (note that for $\kappa=\lbb i\rbb$
we have $\fEC_\alpha^\kappa = \fE_\alpha^i$). To this end we consider
the following subspaces
$$
\begin{array}{l}
\fEF = p_{12}^{-1}(\fE_{\alpha}^i)\cap p_{23}^{-1}((\fE_{\alpha+i-j}^j)^T)
\subset \CMPC\alpha\times\CMPC{\alpha+i}\times\CMPC{\alpha+i-j},
\smallskip\\
\fFE = p_{12}^{-1}((\fE_{\alpha-j}^j)^T)\cap p_{23}^{-1}(\fE_{\alpha-j}^i)
\subset \CMPC\alpha\times\CMPC{\alpha-j}\times\CMPC{\alpha+i-j}.
\end{array}
$$
Consider the following open subset
$$
U = \begin{cases}
\CMPC\alpha\times\CMPC{\alpha+i-j}, & \text{if $i\ne j$}\\
\CMPC\alpha\times\CMPC{\alpha} - \Delta & \text{if $i=j$}
\end{cases}
$$
where $\Delta$ is the diagonal. Then it is easy to see that
$$
\fEF\cap p_{13}^{-1}(U)\cong \fFE\cap p_{13}^{-1}(U).
$$
Indeed, the map
$$
\fFE\cap p_{13}^{-1}(U)\ni(\CF_\bullet\subset\CF'_\bullet\supset\CF''_\bullet)
\mapsto
(\CF_\bullet\supset(\CF_\bullet\cap\CF''_\bullet)\subset\CF''_\bullet)
\in\fEF\cap p_{13}^{-1}(U)
$$
gives such an isomorphism. Hence we have
$$
[e_i,f_j]=0\quad\text{for $i\ne j$,}\quad\text{and}\quad
[e_i,f_i]=b_\alpha^i[\Delta]\quad\text{for $i\ne j$,}
$$
and it remains to compute $b_\alpha^i$.

So, assume that $i=j$. Let us begin with the contribution
of $\fFE$ into $b_\alpha^i$. To this end, note that the fiber
of $\fFE$ over generic point of the diagonal $\Delta$ (with respect
to the projection $p_{13}$) is empty. The reason is the fact
that for a locally free $\CF_\bullet$ there exists no $\CF'_\bullet$,
such that $\CF_\bullet\subset\CF'_\bullet$. Thus $\fFE$ doesn't
contribute into $b_\alpha^i$.

As for $\fEF$, the situation here is quite opposite. For generic point
$\xi = (\CF_\bullet,\CF_\bullet)\in\Delta\subset\CMPC\alpha\times\CMPC\alpha$
the fiber of $\fEF$ over $\xi$ is isomorphic to $\bC\cong\BP^1$: it consists
of all $\CF'_\bullet\subset\CF_\bullet$ such that
$\kappa(\CF_\bullet/\CF'_\bullet)=\lbb i\rbb$, and such subobjects
are uniquely determined by the point
$x=\supp(\CF_\bullet/\CF'_\bullet)=\br(\CF_\bullet,\CF'_\bullet)$.
Moreover, the intersection
$p_{12}^{-1}(\fE_{\alpha}^i)\cap p_{23}^{-1}((\fE_{\alpha+i-j}^j)^T)$
in this case has dimension greater by~1 than expected, thus we
are in the excess intersection situation. It follows that
$b_\alpha^i$ equals to the degree of the excess intersection line bundle
restricted to the fiber $\fEF_\xi$. Further, acting as in \cite{fk1}~3.6.1
we can show that
$$
b_\alpha^i = \deg\bq_\xi^*\CN_{D_\alpha^i/\CMPC{\alpha+i}},
$$
where $D_\alpha^i=\bq(\fE_\alpha^i)$, and 
$\bq_\xi:\bp^{-1}(\CF_\bullet)\cong\bC \subset \fE_\alpha^i 
\to \CMPC{\alpha+i}$ is the canonical projection.

Now let us identify the normal bundle $\CN_{D_\alpha^i/\CMPC{\alpha+i}}$.
Let $\CF_\bullet\in\CMPC{\alpha}$ be a locally free parabolic sheaf and
assume that $\CF'_\bullet\subset\CF_\bullet$ is such that 
$\kappa(\CF_\bullet/\CF'_\bullet)=\lbb i\rbb$. Let  
$c=\supp(\CF_\bullet/\CF'_\bullet)=\br(\CF_\bullet,\CF'_\bullet)$.
Then $\CF'_\bullet\in D_\alpha^i$, and we have the following exact sequence
\begin{equation}\label{cfseq}
0 \to \CF'_\bullet \to \CF_\bullet \to 
(\CF_i/\CF_{i-1})_c\otimes_\BC\CO_c[i] \to 0,
\end{equation}
where the right term is considered as an $n$-periodic representation 
of an infinite linear quiver in the category of sheaves on $\CS'$ with
the sheaves $(\CF_i/\CF_{i-1})_c\otimes_\BC\CO_c$ placed at
$k\equiv i \pmod n$ and with zero at all other places $k$.
Now we want to compute the tangent space (see ~\cite{y})
$$
\CT_{\CF'_\bullet}\CMPC{\alpha+i} = 
\Ext^1(\CF'_\bullet,\CF'_\bullet(-\bD_\infty))
$$
using the exact sequence (\ref{cfseq}). Here $\bD_\infty$ stands for 
$\bC\times\bx$.
To this end we have to compute 
$$
\begin{array}{l}
\Ext^*((\CF_i/\CF_{i-1})_c\otimes_\BC\CO_c[i],\CF_\bullet(-\bD_\infty)),\\
\Ext^*((\CF_i/\CF_{i-1})_c\otimes_\BC\CO_c[i],
(\CF_i/\CF_{i-1})_c\otimes_\BC\CO_c[i](-\bD_\infty)),\\
\Ext^*(\CF_\bullet,(\CF_i/\CF_{i-1})_c\otimes_\BC\CO_c[i](-\bD_\infty)).
\end{array}
$$
The third $\Ext$ is easiest to compute. It is clear that we have
$$
(\CF_i/\CF_{i-1})_c\otimes_\BC\CO_c[i](-\bD_\infty) \cong 
(\CF_i/\CF_{i-1})_c\otimes_\BC\CO_c[i]
$$
and
$$
\Hom(\CF_\bullet,(\CF_i/\CF_{i-1})_c\otimes_\BC\CO_c[i]) =
(\CF_i/\CF_{i-1})_c\otimes_\BC(\CF_i/\CF_{i-1})^*_c \cong \BC,\quad
\Ext^{>0} = 0.
$$
To compute the other two $\Ext$-s we use the following locally free 
resolution of $\CO_c[i]$:
%$$
%\begin{CD}
%@.     @AAA                   @AAA                              @AAA     @AAA \\
%0 @>>> \CO(-\bX_c)       @>>> \CO(-\bX_c)\oplus\CO               @>>> \CO @>>> 0 \\
%@.     @AAA                   @AAA                              @AAA     @AAA \\
%0 @>>> \CO(-\bX_c-\bD_0) @>>> \CO(-\bX_c)\oplus\CO(-\bD_0)       @>>> \CO   @>>> \CO_c \\
%@.     @AAA                   @AAA                              @AAA     @AAA \\
%0 @>>> \CO(-\bX_c-\bD_0) @>>> \CO(-\bX_c-\bD_0)\oplus\CO(-\bD_0) @>>> \CO(-\bD_0) @>>> 0 \\
%@.     @AAA                   @AAA                              @AAA     @AAA 
%\end{CD}      
%$$
$$
\arraycolsep=0pt
\begin{array}{ccccrclcccccc}
&& \uparrow &&& \uparrow &&& \uparrow && \uparrow \\
0 & \ \to\  & \CO(-\bX_c) & \ \to\  & \CO(-\bX_c)&\,\oplus\,&\CO & \ \to\  &
\CO & \ \to\  & 0 & \ \to\  & 0\\
&& \uparrow &&& \uparrow &&& \uparrow && \uparrow \\
0 & \to & \CO(-\bX_c-\bD_0) & \to & \CO(-\bX_c)&\oplus&\CO(-\bD_0) & \to & 
\CO & \to & \CO_c & \ \to\  & 0 \\
&& \uparrow &&& \uparrow &&& \uparrow && \uparrow \\
0 & \to & \CO(-\bX_c-\bD_0) & \to & \CO(-\bX_c-\bD_0)&\oplus&\CO(-\bD_0) & \to
& \CO(-\bD_0) & \to & 0 & \ \to\  & 0  \\
&& \uparrow &&& \uparrow &&& \uparrow && \uparrow 
\end{array}
$$
Here $\bX_c$ stands for $c\times\bX$.
The rows of the above diagram are exact sequences of coherent sheaves 
and the columns are $n$-quasi-periodic representations of an infinite 
linear quiver. The quasi-periodicity means that when one shifts to 
$n$ positions up, the sheaf became twisted by $\CO(\bD_0)$.

Using this resolution one can easily compute $\Ext$-s:
\begin{multline*}
\Ext^p((\CF_i/\CF_{i-1})_c\otimes_\BC\CO_c[i],\CF_\bullet(-\bD_\infty)) = \\
\begin{cases}
(\CF_i/\CF_{i-1})^*_c\otimes_\BC(\CF_{i+1}/\CF_{i})_c\otimes_\BC\CO(\bX_c)_c,
& p=2\\
0, & \text{otherwise}
\end{cases}
\end{multline*}
$$
\Ext^p((\CF_i/\CF_{i-1})_c\otimes_\BC\CO_c[i],
(\CF_i/\CF_{i-1})_c\otimes_\BC\CO_c[i](-\bD_\infty)) = 
\begin{cases}
\BC, & p = 0, 1;\\
0,   & \text{otherwise}
\end{cases}
$$
Now, computing $\Ext^*(\CF'_\bullet,\CF'_\bullet(-\bD_\infty))$
with the help of (\ref{cfseq}) one gets a spectral sequence with 
the first term as follows:
$$
\xymatrix@R=0pt{
{(\CF_i/\CF_{i-1})^*_c\otimes_\BC(\CF_{i+1}/\CF_{i})_c\otimes_\BC\CO(\bX_c)_c} 
& 0 & 0 \\
0 & {\Ext^1(\CF_\bullet,\CF_\bullet(-\bD_\infty))\oplus\BC} & 0 \\
0 & \BC \arrow[r] & \BC
}
$$
Here the map in the bottom row is the map 
\begin{multline*}
\Hom((\CF_i/\CF_{i-1})_c\otimes_\BC\CO_c[i],
(\CF_i/\CF_{i-1})_c\otimes_\BC\CO_c[i](-\bD_\infty)) \to \\
\Hom(\CF_\bullet,(\CF_i/\CF_{i-1})_c\otimes_\BC\CO_c[i](-\bD_\infty))
\end{multline*}
induced by the projection 
$\CF_\bullet \to (\CF_i/\CF_{i-1})_c\otimes_\BC\CO_c[i]$.
It is clear that it takes the identity homomorphism to 
this projection. Thus this map is not trivial, hence induces
an isomorphism in the bottom row of the spectral sequence.
It follows that the spectral sequence degenerates in the second term
and gives the following exact sequence
\begin{multline*}
0 \to \Ext^1(\CF_\bullet,\CF_\bullet(-\bD_\infty))\oplus\BC \to
\Ext^1(\CF'_\bullet,\CF'_\bullet(-\bD_\infty)) \to \\
(\CF_i/\CF_{i-1})^*_c\otimes_\BC(\CF_{i+1}/\CF_{i})_c\otimes_\BC\CO(\bX_c)_c 
\to 0.
\end{multline*}
It is clear that the first term in this exact sequence is the tangent
space to the divisor $D_\alpha^i$ at the point $\CF'_\bullet\in\CMPC\alpha$. 
Hence, the fiber of the normal bundle at this point is isomorphic to
\begin{equation}\label{fibcn}
(\CF_i/\CF_{i-1})^*_c\otimes_\BC(\CF_{i+1}/\CF_{i})_c\otimes_\BC\CO(\bX_c)_c.
\end{equation}

Now we can compute $b_\alpha^i$. To this end we should let the point $c$
vary within the curve $\bC$ and compute the degree of the line bundle
formed by spaces (\ref{fibcn}). The bundle in question is clearly
$$
(\CF_i/\CF_{i-1})^*\otimes(\CF_{i+1}/\CF_{i})\otimes\CO(2),
$$
the last factor is the restriction of the sheaf $\CO(\Delta_\bC)$ 
on $\bC\times\bC$ to the diagonal. Thus
$$
b_\alpha^i = \deg((\CF_i/\CF_{i-1})^*\otimes(\CF_{i+1}/\CF_{i})\otimes\CO(2))
= - \deg(\CF_i/\CF_{i-1}) + \deg(\CF_{i+1}/\CF_{i}) + 2.
$$
Applying Lemma~3.1.1 from \cite{fk2} we get
$$
b_\alpha^i = \langle i',\alpha\rangle + 2. 
%= \langle i',\alpha + 2\rho\rangle.
$$
This completes the proof of the Proposition.
\qed

\medskip

Recall now that $\hsl_n$ is a Kac-Moody algebra with generators $e_i,f_i,h_i,\
i\in I$, and standard relations. 
Theorem ~\ref{hallth} together with Proposition ~\ref{ef} combine into
the following

\begin{theorem}
The operators $e_i,f_i,h_i,\ i\in I$, generate an action of $\hsl_n$ 
on $\bigoplus_\alpha H^\bullet(\CM^\alpha_{\mathbf{gt}})$. This action has
central charge 2.
\end{theorem}

{\sl Proof:} 
It is well known that the subalgebra $U^+\subset U(\hsl_n)$ generated
by $e_i,\ i\in I$, 
embeds into the Hall algebra $\bH,\ e_i\mapsto S_{\{\{i\}\}}$. Thus the 
Serre relations for $e_i, i\in I$, follow.
It only remains to check the Serre relations for the operators $f_i,\ i\in I$. 
But they are given by correspondences transpose
to those of $e_i$.
\qed


\begin{thebibliography}{XXX}

%\raggedright

\bibitem{a} M.F.Atiyah, {\em Instantons in two and four dimensions},
Comm. Math. Phys. {\bf 93}, No. 4 (1984), 437--451.

\bibitem{b} V.Baranovsky, {\em Moduli of Sheaves on Surfaces and Action
of the Oscillator Algebra}, preprint math.AG/9811092 (1998).

\bibitem{bgk} V.~Baranovsky, V.~Ginzburg, A.~Kuznetsov,
{\em Quiver varieties and a noncommutative quadric},
in preparation.

\bibitem{bbe} A.Beilinson, S.Bloch, H.Esnault, {\em $\varepsilon$-factors
for Gauss-Manin determinants}, preprint math.AG/0111277.

\bibitem{bgfm} A.Braverman, M.Finkelberg, D.Gaitsgory, I.Mirkovi\'c,
{\em Intersection Cohomology of Drinfeld's compactifications},
preprint math.AG/0012129.


\bibitem{fg} A.Braverman, M.Finkelberg, D.Gaitsgory,
{\em Crystals and rational curves
in the Kashiwara flag schemes}, in preparation.

\bibitem{fk1} M.Finkelberg, A.Kuznetsov, {\em Global Intersection Cohomology
of Quasimaps' spaces}, Intern. Math. Res. Notices {\bf 7} (1997), 301--328.

\bibitem{fk2} M.Finkelberg, A.Kuznetsov, {\em Parabolic sheaves on surfaces
and affine Lie algebra $\hgl_n$}, J. reine angew. Math. {\bf 529} (2000),
155--203.

\bibitem{fm} E.Frenkel, E.Mukhin, {\em The Hopf algebra
$Rep U_q\widehat{gl}_\infty$}, preprint math.QA/0103126.



\bibitem{g} D.Gaitsgory, {\em Construction of central elements in the affine
Hecke algebra via nearby cycles}, Invent. Math. {\bf 144} (2001),
253--280.

\bibitem{hl} D.Huybrechts, M.Lehn, {\em Stable pairs on curves and
surfaces}, J. Algebraic Geom. {\bf 4} (1995), 67--104.

\bibitem{k} M.Kashiwara, {\em The flag manifold of Kac-Moody Lie algebra},
in Algebraic Analysis, Geometry, and Number Theory, Proceedings of the
JAMI Inaugural Conference, The Johns Hopkins University Press (1989),
161--190.

\bibitem{k2} M.Kashiwara, {\em Kazhdan-Lusztig conjecture for a symmetrizable
Kac-Moody Lie algebra}, Progr. Math. {\bf 87} (1990), 407--433.

\bibitem{ku} A.Kuznetsov, {\em The Laumon's resolution of the Drinfeld's
compactification is small}, Math. Res. Letters {\bf 4}, No. 2--3 (1997),
349--364.

\bibitem{la} G.Laumon, {\em Faisceaux Automorphes Li\'es aux S\'eries
d'Eisenstein}, Perspect. Math. {\bf 10} (1990), 227--281.

\bibitem{l} G.Lusztig, {\em On Quiver Varieties}, Advances in Math. {\bf 136}
(1998), 141--182.

\bibitem{n0} H.Nakajima, {\em Quiver varieties and Kac-Moody algebras},
Duke Math. Journal {\bf 91}, No. 3 (1998), 515--560.

\bibitem{n1} H.Nakajima, {\em Lectures on Hilbert schemes of points on
surfaces}, The AMS University Lecture Series {\bf 18} (1999).

\bibitem{nk} H.Nakajima, {\em Moduli of sheaves on blown-up surfaces},
preprint at\\
http://www.kusm.kyoto-u.ac.jp/$\,\widetilde{}\,$nakajima/TeX.html (2000).

\bibitem{r} C.M.Ringel, {\em The composition algebra of a cyclic quiver},
Proc. London Math. Soc. (3) {\bf 66}, No. 3 (1993), 507--537.


\bibitem{y} K.Yokogawa, {\em Infinitesimal deformation of parabolic Higgs
sheaves}, Intern. J. Math. {\bf 6}, No. 1 (1995), 125--148.

\end{thebibliography}
\end{document}